\documentclass[11pt,reqno]{amsart}
\usepackage[margin=2.5cm]{geometry}
\usepackage{graphicx} 
\usepackage{float}
\usepackage{caption}
\usepackage{subcaption}

\newtheorem{remark}{Remark}[section]

\usepackage{enumitem}
\usepackage{xcolor}

\allowdisplaybreaks

\usepackage{tikz}
\usetikzlibrary{calc}

\usepackage{amsmath}
\usepackage{empheq}
\usepackage{comment}
 \usepackage{amsaddr}
 \usepackage{pifont}
 \newcommand{\cmark}{\ding{51}}%
\newcommand{\xmark}{\ding{55}}%

\usepackage{accents}

\usepackage{caption}
\usepackage{subcaption}
 
\usepackage{appendix}
\usepackage{algorithm}
\usepackage{algpseudocode}

\usepackage{tabularx}
\usepackage{multicol, multirow}
\usepackage{mathtools, stmaryrd}
\usepackage{xparse}
\DeclarePairedDelimiterX{\Iintv}[1]{\llbracket}{\rrbracket}{\iintvargs{#1}}
\NewDocumentCommand{\iintvargs}{>{\SplitArgument{1}{,}}m}
{\iintvargsaux#1} %
\NewDocumentCommand{\iintvargsaux}{mm} {#1\mkern1.5mu,\mkern1.5mu#2}

\usepackage[backend=bibtex,
style = ieee,
]{biblatex}
\usepackage{bbm}
\usepackage{csvsimple}

\def\bSigma{\boldsymbol{\Sigma}}
\def\bTheta{\boldsymbol{\Theta}}

\def\NN{\mathbb N}

\usepackage[normalem]{ulem}
\title{Coincident Peak Prediction for Capacity and Transmission Charge Reduction}

\author{René Carmona, Xinshuo Yang, and Claire Zeng}
\address[A1,A2,A3]{Princeton University, Princeton, USA}
\email[A1,A2,A3]{rcarmona@princeton.edu, xy3134@princeton.edu, cszeng@princeton.edu}
\date{\today}

\begin{document}

\begin{abstract}
    Meeting the ever-growing needs of the power grid requires constant infrastructure enhancement. There are two important aspects for a grid’s ability to ensure continuous and reliable electricity delivery to consumers: capacity, the maximum amount the system can handle and transmission, the infrastructure necessary to deliver electricity across the network. These capacity and transmission costs are then allocated to the end-users according to the cost causation principle. These charges are computed based on the customer's demand on coincident peak (CP) events, time intervals when the system-wide electric load is highest. We tackle the problem of predicting CP events based on actual load and forecast data on the load of different jurisdictions. In particular, we identify two main use cases depending on the availability of a forecast. Our approach generates scenarios and formulates Monte-Carlo estimators for predicting CP-day and exact CP-hour events. Finally, we backtest the prediction performance of strategies with adaptive threshold for the prediction task. This analysis enables us to derive practical implications for load curtailment through Battery Energy Storage System (BESS) solutions. 
\end{abstract}

\maketitle


\section{Introduction, motivations and related literature}

Independent System Operators (ISO) oversee the management of the power grid, in particular the planning of the generation, transmission, and distribution activities. This requires constant enhancement of the infrastructure, including for instance the high-voltage power lines, and other operational expenses that help ensure continuous and reliable delivery of energy to end-consumers. First established by the Federal Energy Regulatory Commission (FERC) in the  Committee Order No. 890,  and later reaffirmed by the DC Circuit, the cost causation principle states that the expenses for new electricity infrastructure must be incurred by the beneficiaries. As a result, capacity and transmission charges are reflected on end-users bills. In particular, commercial and industrial customers face significant peak demand charges that are directly related to their energy usage during time slots when the demand on the grid is the highest. The cost allocation method depends heavily on the ISO or on the utility, but a common form is that of a coincident peak program; a coincident peak (CP) corresponds to the time interval during which the greatest volume of electricity is recorded across the grid according to predefined rules. Notice that the customer's own peak load may not necessarily coincide with any of the system or utility peaks. These methodologies are meant to fairly distribute costs so that it reflects a customer's usage of the transmission infrastructure. 

Understanding capacity and transmission charges is essential for all participants in the power system: indeed, they can represent on average between 10\% and 40\% of a commercial customer's bill. Reducing the demand for electricity during the CP events can hence lead to significant bill savings. Demand Side Management can be achieved through two means: load curtailment and local electricity generation/reserve. For instance, Direct Load Control can involve turning-off non-essential lighting, modifying manufacturing processes, or adjusting the HVAC equipment, while additional generation can involve dialing fossil fuel back pumps with short start and run times, using local Photo-Voltaic (PV) panels or using stored energy. 


\subsection{Overview of Coincident Peak Pricing Programs for Capacity and Transmission Charges}
\vspace{\baselineskip}
We present an overview of some of the coincident peak programs used by ISOs to allocate capacity and transmission costs.  The methods differ in their definition of the base time interval, the time frame, and the jurisdiction (ISO or utility zone) over which the CPs are computed. Some ISOs have a unique coincident peak program to allocate both capacity and transmission charges, such as ERCOT, while others have distinct programs run at two scales: the ISO-level and a more local zone-level. 

\subsubsection{CAISO} 
In the California ISO (CAISO), generation and capacity charges are allocated to grid users through a 12CP program, where each coincident peak is the highest 15 minute interval of the system-wide load usage during the peak hours of a given billing month. CAISO also regroups three Investor Owned Utilities (IOU): Pacific Gas and Electric (PG\&E), San Diego Gas and Electric (SDGE) and Southern California Edison (SCE). There is a specific Time-Of-Use (TOU) classification of the hours of a day distributed between Super-Peak, Peak, Off-Peak and Super Off-Peak according to the ability of each hour to see the CAISO system-wide peak. The Super-Peak hours are specific to only weekdays of July and August. The transmission charges are then distributed to the IOUs, which in turn pass these costs to each customer class (residential, small/medium and large Commercial \& Industrial customer) based on a 12CP methodology. Each CP is computed as the highest 15 min load contribution of the customer's class to the IOU peak over the billing month. The costs are then allocated within each class through either a flat volumetric rate or a non-coincident charge (relying on a customer's own maximum 15-min demand in the month, irrespective of the system-wide condition). 

\subsubsection{ERCOT} 
The Electric Reliability Council of Texas (ERCOT) is the Regional Transmission Organization (RTO) operating Texas's electrical grid. It allocates both capacity and transmission costs according to a four coincident peak (4CP) program for the summer period to compute the capacity charge allocation factor to be applied to the following year. For each month of the summer period (June, July, August and September), the RTO records the system load for each 15 min period of each day, as well as the 15 min interval and day during which the overall month maximum load was achieved; there is therefore a single coincident peak for each month. 

\subsubsection{ISONE}  
The New England ISO (ISONE) oversees the operation of the North-East part of the United States, comprising Connecticut, Maine, Massachusetts, New Hampshire, Rhode Island and Vermont. New England has a 1CP program that relies on the single maximum load demand hour computed on non-holiday weekdays over the whole year spanning from June to May. Historically, ISONE has been summer-peaking with the 1CP occurring only in June, July and August. ISONE allocates transmission costs under the charge labeled Regional Network Service (RNS), which are determined through a 12CP method. Each coincident peak corresponds to a customer's monthly regional network load that is the maximum hour of use coincident with the load of each local network peak during a specific billing month. 

\subsubsection{NYISO} \label{subsec:nyiso_cp_program}

The New York State ISO also has a 1CP program for capacity cost allocation with the following constraints: the maximum peak demand hour is computed over the months of July and August only and must be a non-holiday weekday. This 1CP is then used to calculate the Installed Capacity (ICAP) tag, which is the capacity cost allocation to be billed over the May-April planning year. NYISO then uses a two-step transmission cost allocation similar to that of CAISO and relying on a 12CP zonal allocation per service class based on the monthly coincident peak load of the service class. 

\subsubsection{PJM} \label{subsec:PJM_pseg_cp_program}
The Pennsylvania-New Jersey-Maryland Interconnection (PJM)  is the ISO supervising the operations of the electric grid of thirteen states (in their entirety or only partially) and the District of Columbia. To allocate capacity charges, it runs a five coincident peak (5CP) program where the five CP events are determined using the overall system load during the summer season (June 1st - September 30th). More precisely, they correspond to the five highest non-holiday workweek day peak hours. The Peak Load Contribution (PLC) tag is then multiplied by a capacity charge that takes into account a combination of zone and system-specific scaling factors. This charge is then billed monthly to the customer over the next year. The allocation of Network Transmission Integrated Service (NITS) charges varies by utility zone. For instance, the Public Service Enterprise Group, PSE\&G, uses a 1CP method based on the customer’s load during the summer highest peak of the utility-wide load, while other utilities in New Jersey can levy 5 summer or 12 annual zonal peak hours. These are generally labeled as the Network Service Peak Loads (NSPL). A more detailed description of the CP methodology for the different utilities of PJM is available in Table \ref{tab:PJM_LSE_CP} in Appendix \ref{app:PJM_NITS}. 

\vspace{\baselineskip} 

Battery Energy Storage Systems (BESS) are devices that can store and release electricity through electro-chemical cells. They can be charged through power coming from the main power grid or generated by local generators. BESSs are particularly popular because of the wide variety of possible revenue streams, including energy arbitrage, ancillary services such as frequency regulation, emergency backup power, or peak shaving. In grid-connected microgrids, a BESS can be used to reduce the load during peak hours and thus associated demand charges. In peak shaving, the agent usually tries to reduce their own load peak, on which specific charges are based. For instance, a BESS can be scheduled to charge during off-peak hours and discharge during peak ones to flatten the load (e.g. for data centers). Dispatching the BESS to meet part of the load during the CP hours will similarly reduce the Network Service Peak Load (NSPL) tag, even if the transmission charges depend on the utility system peak load times. A single minor agent will therefore have little to no impact on the considered peak events. The main challenge to CP demand-side management lies in the fact that the coincident peak events are only known at the end of the computation period (e.g. at the end of every month for a 12CP method, or at the end of the summer season in some 1CP programs). Predicting the CP events is therefore essential and has been the approach adopted in the literature on Capacity and Transmission charge avoidance programs. 

\subsection{Related Literature}

Academic publications on the subject are few and far between because of the dominance of for-profit consultancies promoting the adoption of proprietary prediction and management solutions to concerned customers. Most of the existing literature about coincident peak programs develop methods to forecast the system load, whether it be the ISO or the utility one. The load forecasting methods include time-series statistical approaches and classification methods to predict whether a day or a specific hour is going to be a coincident peak event. A popular case study is the Ontario (Canada) market where the capacity charges are known as global adjustment charges and are determined through a 5CP method described in \parencite{Global_Adjustment_Ontario}. \parencite{Jiang_Levman_Golab_Nathwani_2014} proposes a heuristic algorithm to predict whether the next day forecasted maximum load will be one of the 5CP events of Ontario using a 14-day short-term load and weather forecasts to compute the top-5 rank probabilities. The classification as a potential CP event is determined through the probability of that day being in the top 5 values being higher than a static threshold. The management of the BESS is then scheduled through hard coded rules such as discharging between a certain time window on predicted CP days. In this type of studies, two types of errors are of interest: False Positives corresponding to days that are predicted to be one of the 5CP but end up not being one and False Negatives corresponding to days that are predicted to not be among the 5CP but end up being one. A high number of False Positives can tire agents since alerts will be too frequent and lead to unnecessary power curtailment, while False Negatives will represent missed opportunities to manage the load and reduce the capacity charges. \parencite{Jiang_Levman_Golab_Nathwani_2016} examines the economic benefits, the efficiency of this incentive policy on system-wide peak load reduction and analyzes the benefits of downgrading to a lower number CP program. \parencite{Ryu_Makanju_Lasek_An_Cercone_2016} uses a Naive Bayesian algorithm to classify future hours as a CP hour or a non-CP one. The predictors used are discrete time attributes such as hour of the day or day of the week, as well as continuous ones including weather information through temperature or humidity levels and load forecasts. The algorithm can accommodate classification of the next day containing a daily peak, a three-hour window containing a peak or a one-hour window being the peak depending on the highest predicted probabilities of each of the next 24 hour being a peak. For instance, \parencite{liu_prediction_2019} evaluates the performance of several classifiers such as parametric models (e.g. Naive Bayes, Support Vector Machine) and non-parametric models (e.g. Random Forests) or from Deep Learning (e.g. Convolutional Neural Network Long Short-Term Memory). \parencite{Kadri_Mohammadi_Awadallah_2020} propose a model to determine the optimal sizing and scheduling of the BESS  to maximize revenues from ancillary services and on bill savings of the global adjustment charges. 

To the best of our knowledge, one of the few papers treating the transmission charge questions is \parencite{Wu_Ma_Fu_Hou_Rehm_Lu_2022} with an application to the North Carolina market, where the local utility company allocates capacity charges through 12 monthly CP events over one year. They leverage the prediction algorithm developed in \parencite{Fu_Zhou_Ma_Hou_Wu_2022} where ensemble methods are trained to predict the probability that the next day contains the monthly peak hour and the probability for each hour of the next day to be the monthly CP hour. The predictors are day-ahead load and weather forecasts at an hourly frequency.  This algorithm is then used in a two-step scheme rolled daily through each month: \textit{1)} Predict the probability that the next day contains the 1CP event. \textit{2)} If that probability is above a given threshold among $\{0, 2\%, 10\%\}$, optimize the BESS dispatch so as to maximize the average of the hourly net energy exchange of the BESS weighted by the probability of each of the next 24 hours to be the 1CP hour while taking into account characteristics and physical constraints of the BESS. The performance on test data are then used to produce sensitivity analysis with respect to the threshold and the BESS characteristics, as well as to analyze its incurred degradation. Indeed, frequent charging/discharging cycles can significantly impact the performance of a BESS, shorten its lifespan and decrease the Returns on Investment of these technologies. 

\subsection{Contributions}

This work focuses on predicting coincident peaks that are used to allocate capacity and transmission charges. The upshot of our analysis is also its main merit: despite the quality of its performance, it only requires a small amount of publicly available load data. 
We consider three different cases depending on the number of CPs of interest and the availability of forecasts. The simplest setting is a 1CP program where a single coincident peak is computed at the ISO level and both actual and forecast historical data are published by the ISO. The second setting is a 1CP program where there is no forecast directly available for the zone of interest, but forecast data are available for a larger region including the zone in question. Finally, we also consider a higher number of CPs and analyze a specific 5CP program for which forecasts are available. The use of machine learning models to predict coincident peak events relies on the availability and the quality of many datasets of diverse nature and different sources, including for instance accurate load data and historical weather data and forecasts. This is often a deterrent for small-to-medium sized commercial and industrial consumers. In contrast, our work relies solely on two to three publicly available datasets published from a single source, the ISO of interest. The availability of forecasts enables us to model the deviations of the forecasts from the actual loads from historical data, and to build a Monte Carlo scenario generation engine from such a model. The absence of any direct forecast for the electric load of the zone of interest creates an additional challenge. We resolve the issue by relying on forecasts available for a different region and developing a joint stochastic model for both loads. As a result, the Monte-Carlo scenario generation engine for the zone of interest is based on simulations conditioned on the available region load forecast, after estimating from historical data, the joint distribution of the region and zonal actual loads. We then leverage large sets of scenarios by developing a sequential day-ahead approach to the estimation of the probability that the next day is a CP day, and the probabilities of each hour of that day to be the CP hour. The back-testing of different strategies shows the quality of our estimators and the transposability of our algorithm to different jurisdictions. 

\subsection{Outline of the Paper}
Section \ref{sec:frcst_iso} describes the scenario generation scheme for the electric load of a single region for which forecasts are available. Section \ref{sec:frcst_utility} describes the building blocks of our conditional scenario generation scheme. In each section, we estimate probabilities that the next day is going to be a coincident peak day, and in that case, what are the probabilities of each hour of the day is the coincident peak hour. Section \ref{sec:backtest_implications} defines several strategies and analyze their performance for the coincident peak classification task. Several practical implications are also derived for real-life application of load curtailment strategies. Finally, Section \ref{sec:conclu} summarizes our findings and concludes. 


\section{Scenario Generation for the RTO/ISO Load} \label{sec:frcst_iso}

The goal of this section is to present a stochastic model for the hourly electric load of a unique region (e.g. the RTO or the ISO) for which actual loads and forecasts are available. To that end, we leverage the analysis from \parencite{Carmona_Yang_2024} and adapt it to our setting. The underlying idea behind this procedure is that the forecasts published by the ISO contain valuable and complex information and the modeling of load deviations incorporate it in some way into the Monte Carlo scenarios. 

\subsection{Input of the model}

We assume the availability of the following datasets: 
\begin{itemize}
    \item the actual load historical data of the RTO at an hourly resolution,
    \item the day-ahead forecast data of the RTO at an hourly resolution.
\end{itemize}
To simplify the presentation of this use case, we consider a single forecast with a forecasting horizon that covers all 24 hours of the day the predictions are generated for. 

\subsection{Methodology}

Similarly to \parencite{Carmona_Yang_2024}, we develop a simulation engine, \texttt{PLProb}, to perform Monte Carlo scenario simulations. We describe below the methodology used to generate scenarios from historical data and forecasts. 

Let $\overline{L}^{d}_{h}$ be the actual load at hour $h$ of day $d$ and $\tilde{L}^{d}_{h}$ be the corresponding forecast. Given the chosen historical forecast series, for each day $d$, we read a vector of $N_h = 24$ load points from the actual load and the corresponding load forecasts to create a time series of load deviations defined by: 
\begin{equation}
    L^d_{h} = \overline{L}^{d}_{h} - \tilde{L}^{d}_{h}, \qquad h \in \{0, 1, \dots, N_{h}-1 \}
\end{equation}

Let us denote by $d^*$ the day for which we want to generate scenarios. We denote by $N_d$ the number of days in the history strictly prior to $d^*$. 

\begin{center}

\captionof{algorithm}{Load Deviation Marginal Distribution}\label{algo:load_dev}
\begin{enumerate}[leftmargin=*]
    \item For each time horizon $h  \in \{0, 1, \dots, N_{h}-1 \}$, we have a time series of load deviations $(L^{1}_{h}, L^{2}_{h}, \dots, L^{d}_{h}, \dots, L^{d^*-1}_{h})$ of length $N_d$ and indexed by $d$ such that $d < d^*$, the day the scenarios are generated for. We create a data matrix with $N_d$ rows and $N_h$ columns where the coefficient of position $(i,j)$ is the load deviation at hour $j$ of day $i$. 
    \item For each time horizon $h$, a Generalized Pareto Distribution (GPD) is fitted to the load deviations. We denote such a distribution as $G_{h}$ for the sake of later reference. Throughout this work, in order to manipulate GPD distributions, we use functions from a \texttt{Python} package wrapping the \texttt{R} library \texttt{Rsafd} \parencite{Rsafd} accompanying the book \parencite{Carmona_2014}. 

    \item For each time horizon $h$, the load deviation time series is then transformed into a uniform time series $\{L^d_{h}\}_{d=1}^{N_h}$ by the probability integral transform. Let the cumulative distribution function (CDF) of the GPD be $\Phi_{G_{h}}$ for every day $d<d^*$, we compute $\Phi_{G_{h}}(L^d_{h})$.  
    
    \item For each time horizon $h$, the uniform time series $\Phi_{G_{h}}(L^d_{h})$ is then transformed into a standard Gaussian marginal distribution $\mathcal{N}(0,1)$ by inversion of the probability integral transform. For
each $d<d^*$, we compute $\hat{L}^d_{h} = Q^{-1}(\Phi_{G_{h}}(L^d_{h}))$ where  $Q^{-1}$ is the quantile function of the standard Gaussian
distribution. 

    \item We fit a Gaussian graphical model to estimate the temporal precision matrix based on the data matrix of dimensions $N_d \times N_h$. The previous steps ensure that each column actually follows a standard Normal distribution and we assume that the $N_h$-dimensional random variables are \textit{jointly Gaussian} with mean zero and a certain covariance matrix $\bSigma$ of dimension $N_h \times N_h$.  Hence $\hat{L}^d := (\hat{L}^d_0, \hat{L}^d_1, \dots, \hat{L}^d_{23})^{\dagger} \sim \mathcal{N}(\mathbf{0}_{N_h}, \bSigma)$. 

    The estimation of the precision matrix is done through an $L_1$ regularised Maximum Likelihood Estimator, also known as the LASSO Gaussian Graphical Model. The precision matrix $\bTheta = \bSigma^{-1}$ is estimated by: 
    \begin{equation}
        \hat{\bTheta} := \underset{\bTheta \in \mathfrak{S}_{N_h}}{\operatorname{argmin}}  \Big\{ \operatorname{trace}(S \bTheta) - \log\det (\bTheta) + \lambda \lVert \bTheta \rVert_1 \Big\} 
    \end{equation}
    where $\mathfrak{S}_{N_h}$ denotes the set of symmetric matrices of dimension $N_h \times N_h$, $S$ is the empirical covariance matrix, $\lVert \cdot \rVert_1$ the $L_1$ norm and $\lambda$ the regularization parameter for the $L_1$-penalty.  Given a dimension $n \in \mathbb{N}^*$, the negative log-determinant function operator is defined for $A \in \mathfrak{S}_n$ as: 
    \begin{equation*}
        - \log\det (A) := \begin{cases}
            - \sum_{k=1}^n \log \mu_k(A) & \text{ if } A \succ 0, \\
            + \infty & \text{ otherwise. } 
        \end{cases}
    \end{equation*}
    where $\{\mu_k(A)\}_{k=1}^n$ is the set of eigenvalues of $A$. 
    \item The final step is to generate Monte Carlo samples of the load. Let $K$ be the size of the scenario batch and $s$ be the index of a single scenario. Given samples $\hat{l}^s_h$, $s \in \{1, \dots, K\}$ drawn from the Gaussian distribution $\mathcal{N}(\mathbf{0}_{N_h}, \hat{\bSigma})$, we can generate samples of the load deviations as $l^s_h = \Phi^{-1}_{G_h} (Q(\hat{l}^s_h))$. Realizations of the system's actual load behavior are obtained by simply adding back these deviations to the original forecast.    

\end{enumerate}
\end{center}

\begin{figure}[htb!]
\includegraphics[width=5cm]{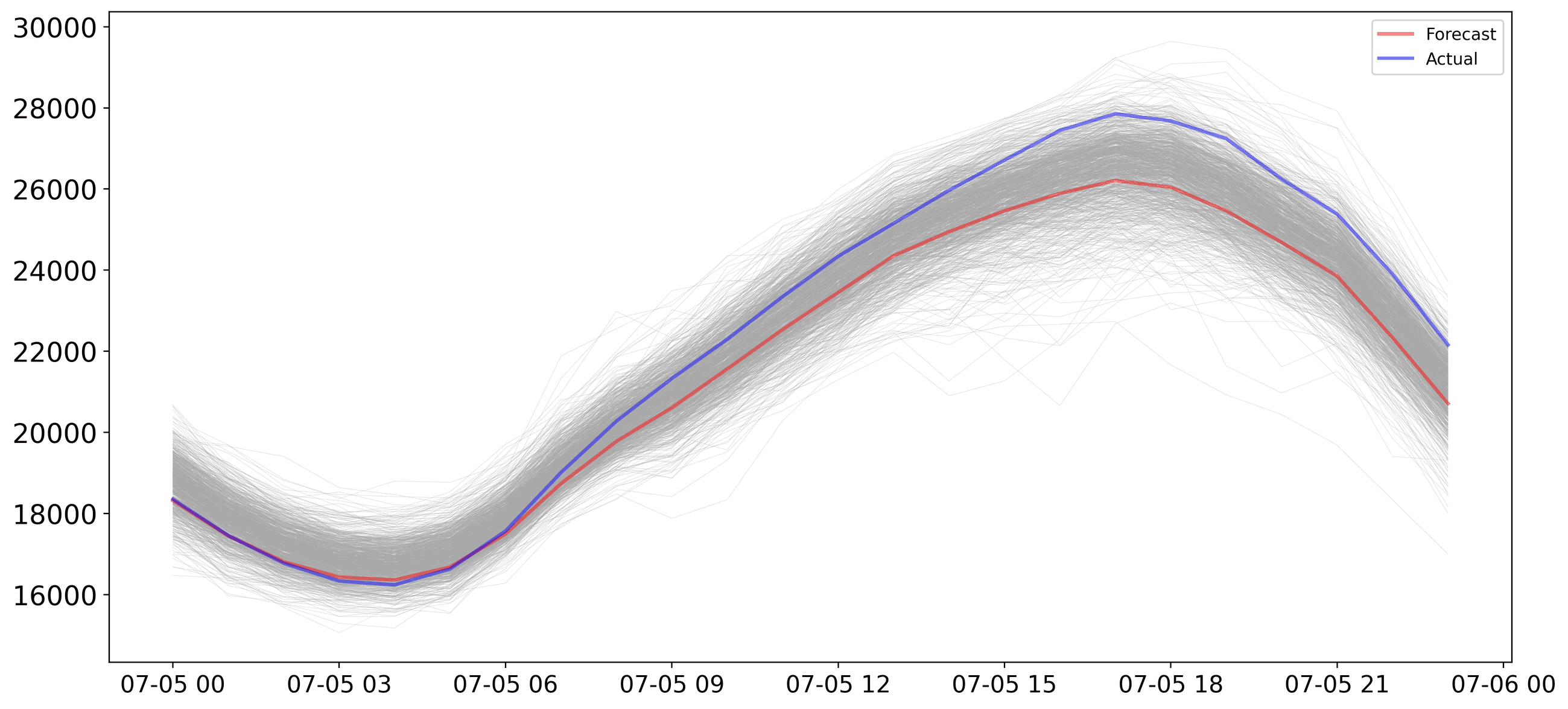}
\includegraphics[width=5cm]{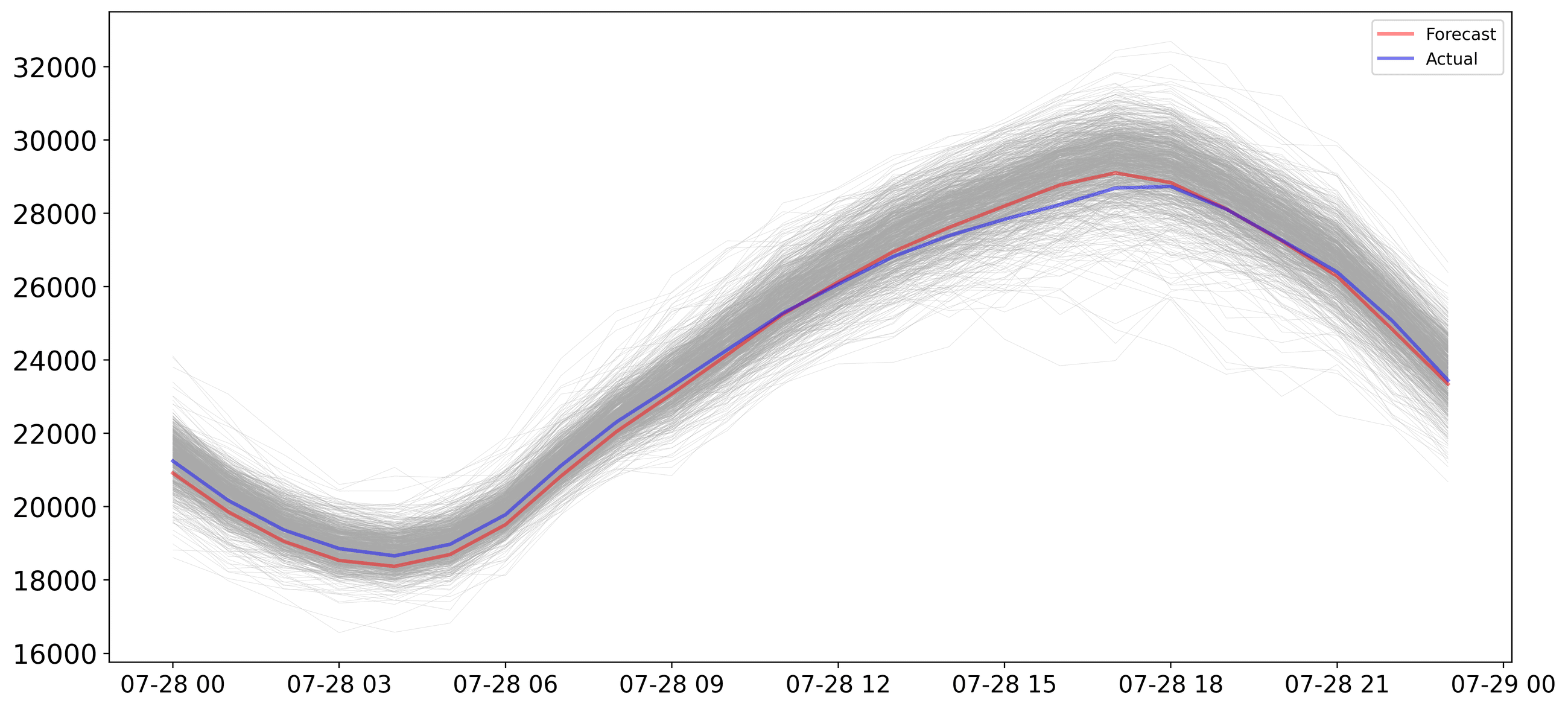}
\includegraphics[width=5cm]{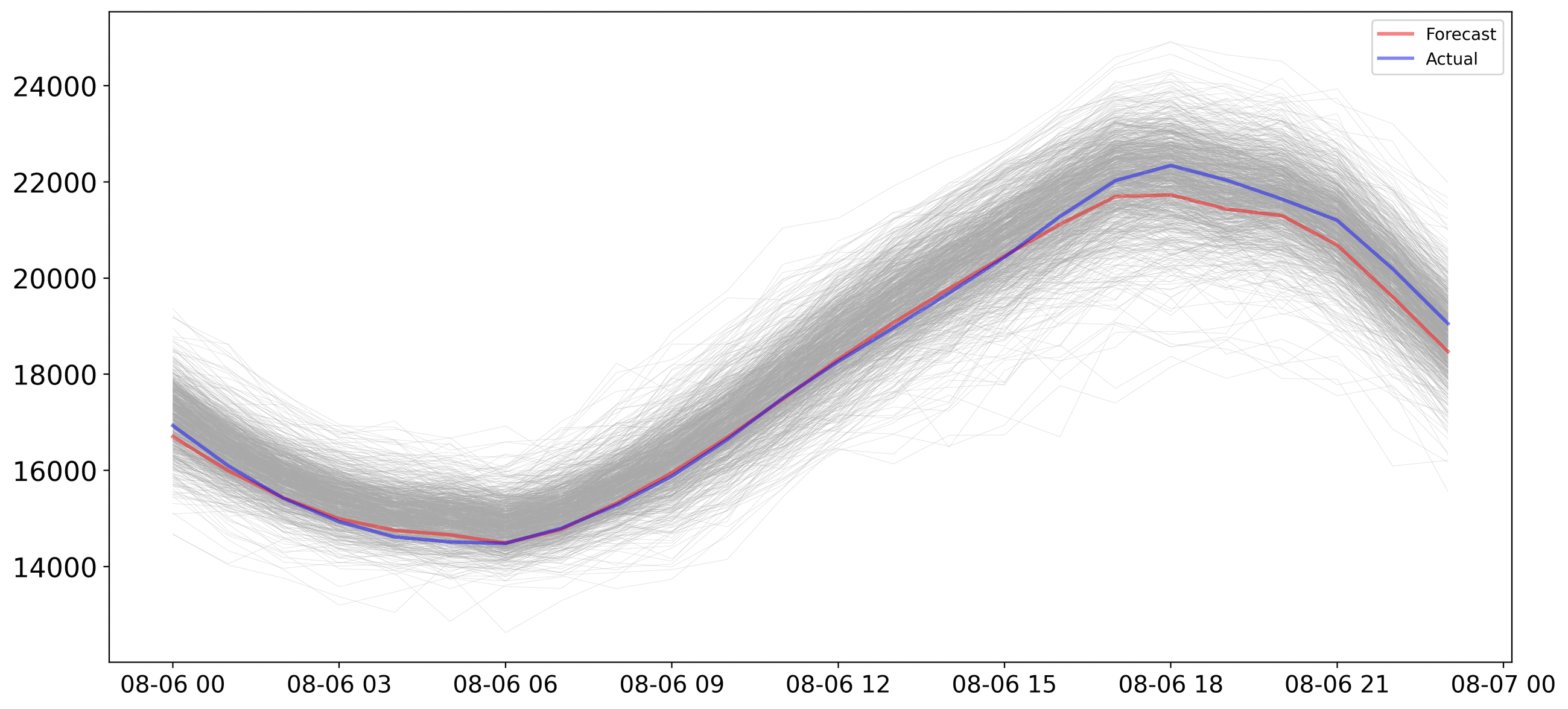}
\caption{PLProb scenarios of NYISO load for different days}
\label{fig:nyiso_scenario}
\end{figure}

To illustrate our scheme, we focus on the NYISO electric load and its 1CP program described in \ref{subsec:nyiso_cp_program}. Figure \ref{fig:nyiso_scenario} illustrates an example of scenarios obtained for the electric load of NYISO on three days with a forecast usually published at 12:00PM (noon) the day before and a forecasting horizon of 24 hours (from midnight till the end of the day). 

\subsection{\textbf{Running maximum probability}}

We now look at the problem of estimating the probability that today is the new running CP day for the current year (in the sense that today's maximum load will be the new running maximum for the system-wide load time-series up until yesterday). 

Let us denote by $d^*$ the current date, by $\{\bar{L}_{d,h}\}_{d \in \NN^*, h \in \{0, \dots, 23\}}$ the RTO actual load on day $d$ during the hour $h$ and by $\{M_{d}\}_{d \in \NN^*}$ the RTO daily maximum actual load of day $d$. We denote by $\operatorname{CP}_{d^*} = \underset{d < d^*}{\operatorname{argmax}} \{ M_d \}$ the running maximum seen up until the day prior to $d^*$. Let $K$ be the number of generated scenarios and we denote by an upper-script $l^k_{d,h}$ the load for day $d$ and hour $h$ and $m^k_d$ maximum daily load of scenario $k \in \{1, \dots, K\}$ of day $d$.

We compute the time-series of probabilities that today's peak is higher than the previous peaks by computing the frequency of scenarios whose daily maximum load exceeds the running maximum of days prior to $d^*$: 
$$ \widehat{\operatorname{prob}}_1 = \frac{1}{K} \sum_{k=1}^K \mathbbm{1}_{\{m^k_{d^*} > \operatorname{CP}_{d^*}\}} $$
with the convention that $\operatorname{CP}_{0} = 0$. 

For the most recent forecast with a forecasting horizon of 12 hours, we plot:  
\begin{itemize}
    \item The actual daily RTO system peak load (red continuous curve)
    \item Green dots indicating the days when the daily maximum load exceeded the running maximum up until now;
    \item For each day, a blue vertical bar whose height is the probability that today's maximum load will exceed the previous ones.
\end{itemize}

\begin{figure}[htb!]
    \centering
    \includegraphics[width = 14cm]{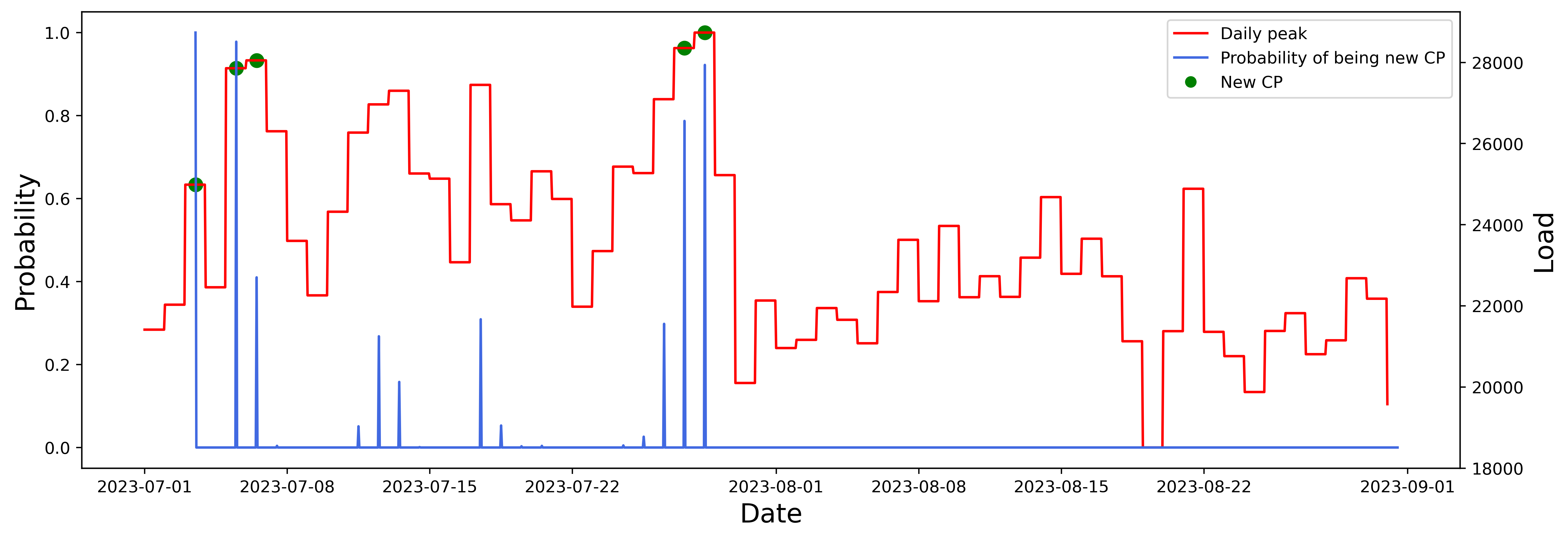}
    \caption{1CP probability for NYISO}
    \label{fig:NYISO_1CP_ts}
\end{figure}

The blue spikes correspond to elevated estimated probabilities that the day in question is a CP day. The first two days of the timeseries are ignored because they correspond to weekend days. We note the presence of probability spikes for all of the coincident peak updates. It is interesting to note that some spikes are less strong (around 0.5), especially if they follow a CP update that ends up having a daily maximum load close to the one they see. A case in point is the relatively low probability assigned to the third update around July 7th in NYISO. This is mainly due to the way the probabilities are computed: the probability of today's load to be a new CP is computed as its probability of surpassing the current running CP. This leads to lower estimates, although any values forecasted to be close to the running CP should be considered due to the stochasticity of the load. These considerations will be taken into account in Section \ref{sec:backtest_implications} to design coincident peak signals. 

\subsection{Estimation of the time of the daily peak}

Given that a day is predicted to be a CP-day, we now need to estimate the probability of each hour being a CP-hour. We are looking for the probability that each hour $h$ is the maximum daily load from the scenarios $k \in \{1, \dots, K\}$ on that day and we denote it by $\operatorname{prob}^{h}_1$. This leads to a vector of length $24$ containing the estimators: 
$$ \widehat{\operatorname{prob}}^{h}_1 = \frac{1}{K} \sum_{k=1}^K \mathbbm{1}_{\{l^k_{d^*, h} = m^k_{d^*}\}}, \quad h \in \{0, \dots, 23\}$$

For a set of days on which the CP is updated during Summer 2023, we plot a $2\times 2$ matrix of panes showing
\begin{itemize}
\item the hourly actual load of the day as a continuous red curve,
\item the histogram of the hourly probabilities that a given hour is the hour the daily maximum occurs, computed from the generated scenarios,
\item a point in green indicating the hour when the maximum load actually occurs. 
\end{itemize}

\begin{figure}[htb!]
    \centering
    \includegraphics[width = 14cm]{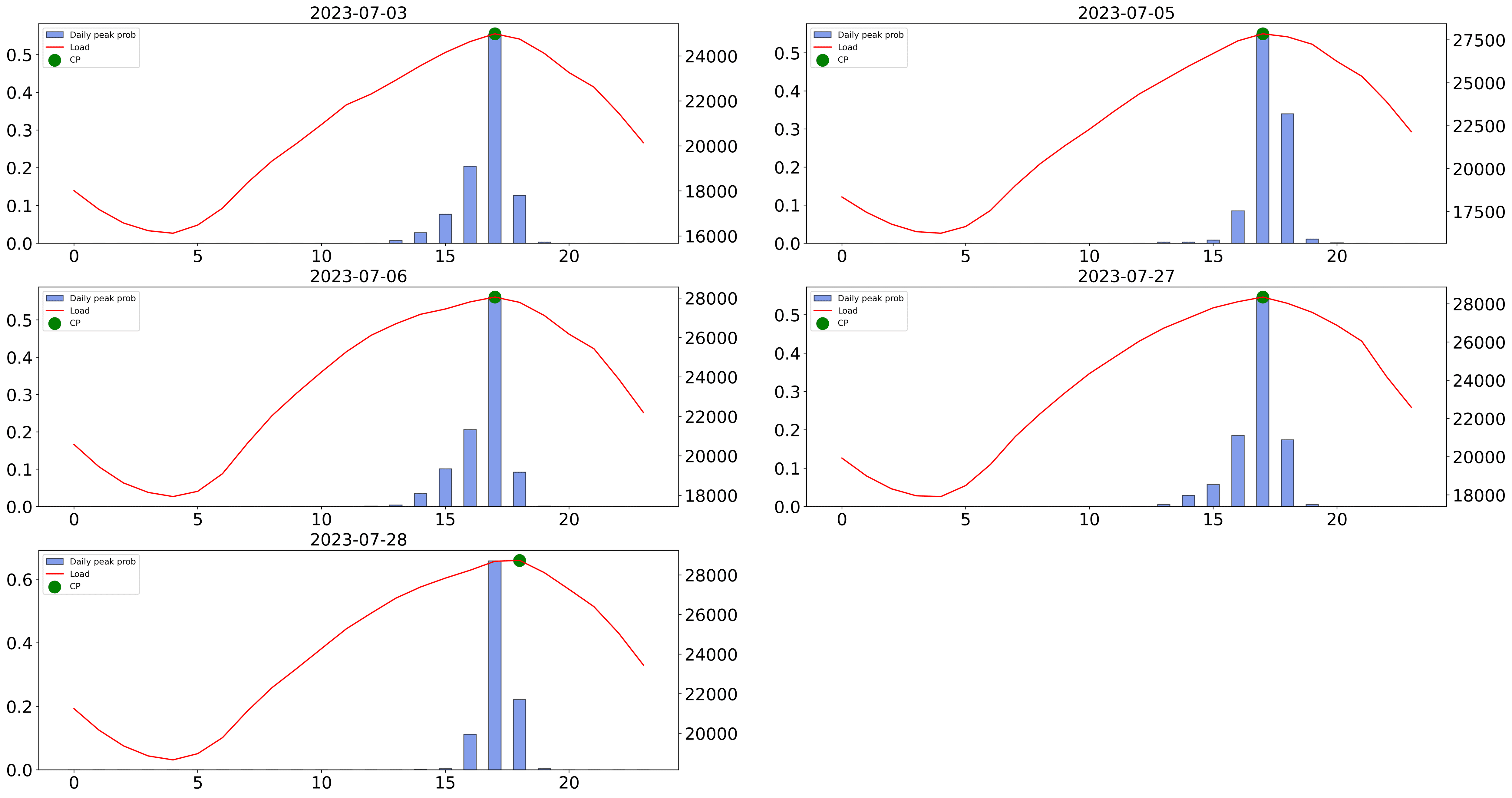}
    \caption{1CP probability for NYISO}
    \label{fig:NYISO_1CP_max_hour}
\end{figure}

We first note that all of the non-zero probabilities are concentrated in the afternoon, which is consistent with what is actually observed afterwards. The hour with the highest assigned probability is almost always among the true top hourly actual loads of the day. Thus this method leads to approximating the time of the CP hour relatively well, with an error of at most 1 hour on the true CP updates.


\section{Conditional Scenario Generation for the Utility Load} \label{sec:frcst_utility}

The goal of this section is to develop a joint stochastic model for the hourly electric load in two different regional zones capturing their correlation, and to design an estimation procedure and a Monte Carlo simulation algorithm capable of producing daily and hourly predictions for the existence of daily and intraday hourly peaks. Motivated by the analysis presented in \parencite{Carmona_Yang_2024}, we propose an alternative modeling and estimation procedure whose implementation does not require the existence of forecasts for all the zonal load time series.


\subsection{Input of the Model} \label{subsec:data_description}
This setting is focused around the electric load of two zones of interest $z_1$ and $z_2$ with the characteristic that zone $z_2$ is a subset of zone $z_1$. We assume that we have the following information: 
\begin{itemize}
    \item the actual load historical data of zone $z_1$ at an hourly resolution,
    \item the actual load historical data of zone $z_2$ at an hourly resolution,
    \item the day-ahead forecast data of zone $z_1$  at an hourly resolution.
\end{itemize}
This structure is summarized in Table \ref{tab:data_pseg}. 

\begin{table}[htb!]
    \centering
    \begin{tabular}{|c|c|c|}
    \hline 
       Zone  &  Actual Load is available & Forecast is available \\
       \hline 
       $z_1$ &  \cmark  & \cmark \\
       $z_2$ & \cmark & \xmark  \\ 
       \hline 
    \end{tabular}
    \caption{Data availability}
    \label{tab:data_pseg}
\end{table}

For each day $d$ and each load zone $z \in \mathcal{Z} = \{z_1, z_2\}$, we have a vector of $N_h$ load point forecasts, with $h \in \{0, 1, \dots, N_h -1\}$ denotes the hour index of the forecast. Let $\overline{L}^{d}_{z,h}$ be the actual load of zone $z$ at hour $h$ of day $d$ and $\tilde{L}^{d}_{z,h}$ be the corresponding forecast. We also consider the load deviations of zone $z_1$: 
\begin{equation}
    L^d_{1,h} = \overline{L}^{d}_{1,h} - \tilde{L}^{d}_{1,h}, \qquad h \in \{0, 1, \dots, N_{h}-1 \}
\end{equation}

\subsection{Conditional Scenario Generation}

We develop a simulation engine \texttt{PLProb} to perform Monte Carlo scenario simulations of the zone $z_2$ load forecast conditioned on the zone $z_1$ forecast. Several building blocks are necessary for this scheme and we describe them below. 

\subsubsection{Marginal Distribution Fitting for each of the actual loads of zones $z_1$ and $z_2$} 
\label{subsubsec:marginal_load}
We fit the marginal distributions of the actual loads of zone $z_1$ and that of zone $z_2$. We present below a simplified version, similar to Algorithm \ref{algo:load_dev} from Section \ref{sec:frcst_iso}. Again, for the purpose of illustration, we assume that we use a forecast with a 24-hour horizon and for each day $d$, we read a vector of $N_h = 24$ load points from the actual load.  Let us denote by $d^*$ the day for which we want to generate scenarios and $z$ the zone we consider. 

\begin{center}
    \captionof{algorithm}{Marginal Distribution Fitting} \label{algo:marginal_fitting}
\begin{enumerate}
    \item For each time horizon $h  \in \{0, 1, \dots,  23 \}$, we have a time series of the hourly loads $\{ \overline{L}^d_{z,h} \}_{d=1}^{N_d}$ with a length $N_d$ and where the data points are indexed by $d$ such that $d < d^*$, the day the scenarios are generated for. 
    \item For each time horizon $h$, a Generalized Pareto Distribution (GPD) is fitted to the loads. We denote such a distribution as $G_{z,h}$ for the sake of later reference. 
    \item For each time horizon $h$, the actual load time series is transformed into a uniform time series e $\Phi_{G_{z,h}}(\overline{L}^d_{z,h})$ by the probability integral transform where the cumulative distribution function (CDF) of the GPD is denoted $\Phi_{G_{z,h}}$ for every day $d<d^*$.  
    \item For each time horizon $h$, the uniform time series $\Phi_{G_{z,h}}(\overline{L}^d_{z,h})$ is then transformed into a standard Gaussian marginal distribution $\mathcal{N}(0,1)$ by inversion of the probability integral transform. For
each $d<d^*$, we compute $\hat{\bar{L}}^d_{z,h} = Q^{-1}(\Phi_{G_{z,h}}(\overline{L}^d_{z,h}))$ where  $Q^{-1}$ is the quantile function of the standard Gaussian distribution. 
\end{enumerate}
\end{center}

\subsubsection{Marginal Distribution Fitting for the load deviations of zone $z_1$}  \label{subsubsec:marginal_deviation}
The same scheme can be applied to the load deviations of zone $z_1$ by replacing $\overline{L}^d_{z,h}$ with the corresponding load forecasting error $L^d_{1,h}$. This would be equivalent to applying Algorithm \ref{algo:load_dev} to the actual and forecast data of zone $z_1$ only. 

\subsubsection{Conditional Scenario Scheme}

We can now describe the full scheme that generates scenarios of the zone $z_2$ load forecast conditioned on the  zone $z_1$ forecast. 

\begin{center}
    \captionof{algorithm}{Conditional Scenario Generation}
\begin{enumerate}
    \item Fit a $N_h = 24$ marginal distribution of the actual load for zone $z_1$ thanks to Algorithm \ref{algo:marginal_fitting} and transform the historical data so that they have a standard Gaussian distribution. We denote this vector by $\hat{\bar{L}}^{d}_{1} := (\hat{\bar{L}}^{d}_{1,0}, \dots, \hat{\bar{L}}^{d}_{1,23})^{\dagger}$;
    \item Do the same for the actual load of $z_2$. We denote this vector by $\hat{\bar{L}}^{d}_{2} := (\hat{\bar{L}}^{d}_{2,0}, \dots, \hat{\bar{L}}^{d}_{2,23})^{\dagger}$.
    We assume that $\hat{L}_1$ and $\hat{L}_2$ are jointly Gaussian, so that \begin{align*}
               \begin{pmatrix}
                    \hat{\bar{L}}_1 \\
                    \hat{\bar{L}}_2
                \end{pmatrix}  \sim \mathcal{N}(0, \mathbf{\Sigma}) \\
    \end{align*}
    
    \item We fit a Gaussian Graphical LASSO model to estimate the spatial and temporal covariance structure of  $z_1$ and  $z_2$ actual daily load vectors. The covariance matrix $\bSigma$ is of dimension $48 \times 48$ with a space dimension of 2 and time dimension of 24 for each zone. 
    $$  \bSigma = \begin{pmatrix}
                \bSigma_{11} & \bSigma_{12} \\
                \bSigma_{21} & \bSigma_{22} \\
                \end{pmatrix}  $$
    with $\bSigma_{ij} \in \mathbb{R}_{24\times 24}$ denoting the covariance between vectors $\hat{L}_i$ and $\hat{L}_j$ for $i,j = 1,2$;
    
    \item Generate scenarios for zone $z_1$ using Subsection \ref{subsubsec:marginal_deviation} by fitting the marginal distribution of deviations $L^d_{1,h}$ and adding them back to the available forecast of  $z_1$;
    \item Transform the  $z_1$ scenarios into Gaussian by using the marginals of (1);
    \item Use the joint Gaussian model of (3) to generate Gaussian load scenarios for  $z_2$ conditionally on the knowledge of the scenarios of (5);
    \item Use the marginal distribution of (2) to get ``real'' scenarios of  $z_2$ with the correct marginals.
\end{enumerate}
\end{center}

\begin{figure}[htb!]
   \begin{tikzpicture}
        \begin{scope}[xshift=0cm]
   \node [rectangle, draw, right=0cm, text width=8.2cm] (eq1) {
        \begin{minipage}{\textwidth}
        \centering \textbf{Actual Load}
            \begin{align*}
                 & (1) \, \textbf{Zone $z_1$} & \bar{L}_1 = {\small \begin{pmatrix}
                    \bar{L}_{1,0} \\
                    \bar{L}_{1,1} \\
                    \cdots \\
                    \bar{L}_{1,23} 
                \end{pmatrix}} \longrightarrow  \hat{\bar{L}}_1 = {\small \begin{pmatrix}
                    \hat{\bar{L}}_{1,0} \\
                    \hat{\bar{L}}_{1,1} \\
                    \cdots \\
                   \hat{\bar{L}}_{1,23} 
                \end{pmatrix}} \\
                & &   \\
                & & {\scriptstyle \bar{L}_{z,h} \mapsto Q^{-1}\circ \, \Phi_{\bar{G}_{z,h}}(\bar{L}_{z,h})} \qquad \\ 
                & &  \\
                &(2) \, \textbf{Zone $z_2$} &  \bar{L}_2 = {\small \begin{pmatrix}
                    \bar{L}_{2,0} \\
                    \bar{L}_{2,1} \\
                    \cdots \\
                    \bar{L}_{2,23} 
                \end{pmatrix}} \longrightarrow \hat{\bar{L}}_2 = {\small \begin{pmatrix}
                    \hat{\bar{L}}_{2,0} \\
                    \hat{\bar{L}}_{2,1} \\
                    \cdots \\
                   \hat{\bar{L}}_{2,23}
                \end{pmatrix}   }
            \end{align*}
        \end{minipage}
    }; 
    \node [rectangle, draw, right=9.5cm, text width=4.5cm] (eq1) {
        \begin{minipage}{\textwidth}
            \centering \textbf{Assumption}

            Jointly Gaussian 

            \begin{align*}
               \begin{pmatrix}
                    \hat{\bar{L}}_1 \\
                    \hat{\bar{L}}_2
                \end{pmatrix}  \sim \mathcal{N}(0, \mathbf{\Sigma}) \\
            \end{align*}

            $(3)$ \texttt{GLASSO} to estimate $\mathbf{\hat{\Sigma}}$

        \end{minipage}
    };

     \node [rectangle, draw, below = 6.5cm, right=0cm, text width=8cm] (eq1) {
        \begin{minipage}{\textwidth}
                \centering \textbf{Zone $z_1$ Gaussianized Load Deviation} 
                
                ${L}_1 = \bar{L}_1 - \tilde{L}_1 $
            \begin{align*}
            (4) \quad  & {L}_1\xrightarrow{ L_{1,h} \mapsto Q^{-1}\circ \, \Phi_{G_{1,h}}(L_{1,h})}  \hat{L}_1 \sim \mathcal{N}(0, \mathbf{S})
            \end{align*}
            \texttt{GLASSO} to estimate $\mathbf{\hat{S}}$

            \medskip
            
            \textbf{Zone $z_1$ Load Deviation Scenario}
            \begin{align*}
            (5) \quad   & {L}^s_1 \sim \mathcal{N}(0, \mathbf{\hat{\mathbf{S}}}) \xrightarrow{ L^s_{1,h} \mapsto \Phi^{-1}_{G_{1,h}}\circ \,  Q (L^s_{1,h})}  \hat{L}^s_1 
            \end{align*}

            Obtain realization through $\bar{l}^s_1 = \hat{L}^s_1 + \tilde{L}_1 $
        \end{minipage}
    };

     \node [rectangle, draw, below = 7.25cm, right=9.25cm, text width=5.3cm] (eq1) {
        \begin{minipage}{\textwidth}
            \centering 
             \textbf{Zone $z_2$ Gaussianized Load Conditional Scenario}
        \begin{align*}
            (6) \quad & L^s_2 | \hat{L}_1 = l^s_1 \sim \mathcal{N}(\mu_{2|1},  \bSigma_{2|1})\\
            & \mu_{2|1} := \bSigma_{12} \bSigma_{22}^{-1} \\
            & \bSigma_{2|1} := \bSigma_{11} -  \bSigma_{12} \bSigma_{22}^{-1}  \bSigma_{21}\\
            & \bSigma = \begin{pmatrix}
                \bSigma_{11} & \bSigma_{12} \\
                \bSigma_{21} & \bSigma_{22} \\
            \end{pmatrix} \\
            & \bSigma_{ij} \in \mathbb{R}_{24 \times 24}, \, i,j = 1,2
        \end{align*}
        
        \bigskip 
        
            \textbf{Zone $z_2$ Load Scenario}
        \begin{align*}
            (7) \quad & l^s_2 \sim \mathcal{N}(\mu_{2|1},  \bSigma_{2|1}) \\ 
            & l^2_s  \xrightarrow{ l^s_{2,h} \mapsto \Phi^{-1}_{G_{2,h}} \circ \, Q (l^s_{2,h})}  \hat{L}^s_2  
        \end{align*}
        \end{minipage}
    };

    \draw (11.75cm,-3.6) -- (11.75cm,-2.0);
    \draw (8.9cm,-2.85) -- (11.75cm,-2.85);
    \draw (8.9cm,-2.85) -- (8.9cm,-6.725);
    \draw (8.25cm,-6.725) -- (8.9cm,-6.725);

        \end{scope}
    \end{tikzpicture}
    \caption{Conditional scenario generation scheme for a forecasting horizon of 24 hours}
    \label{fig:diagram_scenario_gen}
    \end{figure}
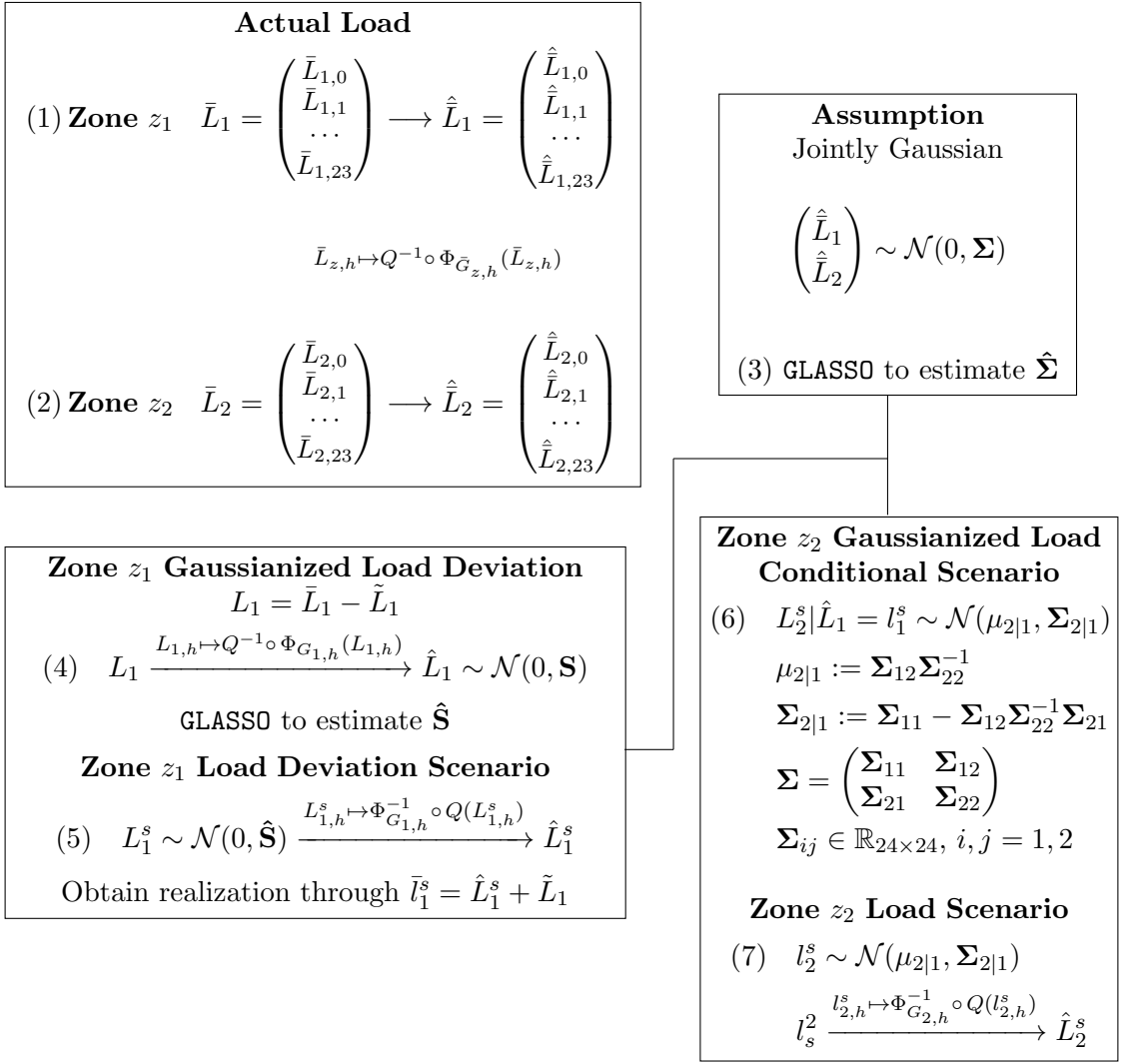 


\subsection{Conditional Simulation Results} 
\label{sec:Conditional Simulation Results}
We illustrate the performance of the model and the Monte Carlo simulation engine presented above by the results of its implementation on publicly available data from PJM \cite{PJM_2011} for the zones $z_1=MIDATL$ and $z_2=PSEG$, the choice of these two particular zones being motivated by the investigation of 1CP program presented in Subsection \ref{subsec:PJM_pseg_cp_program}.

PSE\&G (denoted by PSEG) is part of the Mid-Atlantic (denoted by MIDATL) group of PJM along with utilities described in Table \ref{tab:PJM_LSE_CP}. PJM publishes several datasets of interest on \parencite{PJM_2011}, including: 
\begin{itemize}
    \item the actual load historical data of the Mid-Atlantic region at an hourly resolution,
    \item the actual load historical data of PSE\&G at an hourly resolution,
        \item the day-ahead forecast data of the Mid-Atlantic region  at an hourly resolution.

\end{itemize}
In particular, four different sets of forecasts are published for Mid-Atlantic each day, each with a forecasting horizon that extends beyond 24 hours. We will truncate these forecasts to end at the last hour of the day (11:00 pm).  
\begin{itemize}
\item at 11:45 pm the day before for hours $0,1,\cdots,23$ of the following day,
\item at 5:45 am the day of, for hours $6,7,\cdots,23$ of the day,
\item at 11:45 am the day of, for hours $12,13,\cdots,23$ of the day,
\item at 5:45 pm the day of, for hours $18,19,\cdots,23$ of the day.
\end{itemize}

\begin{figure}[htb!]
\includegraphics[width=5cm]{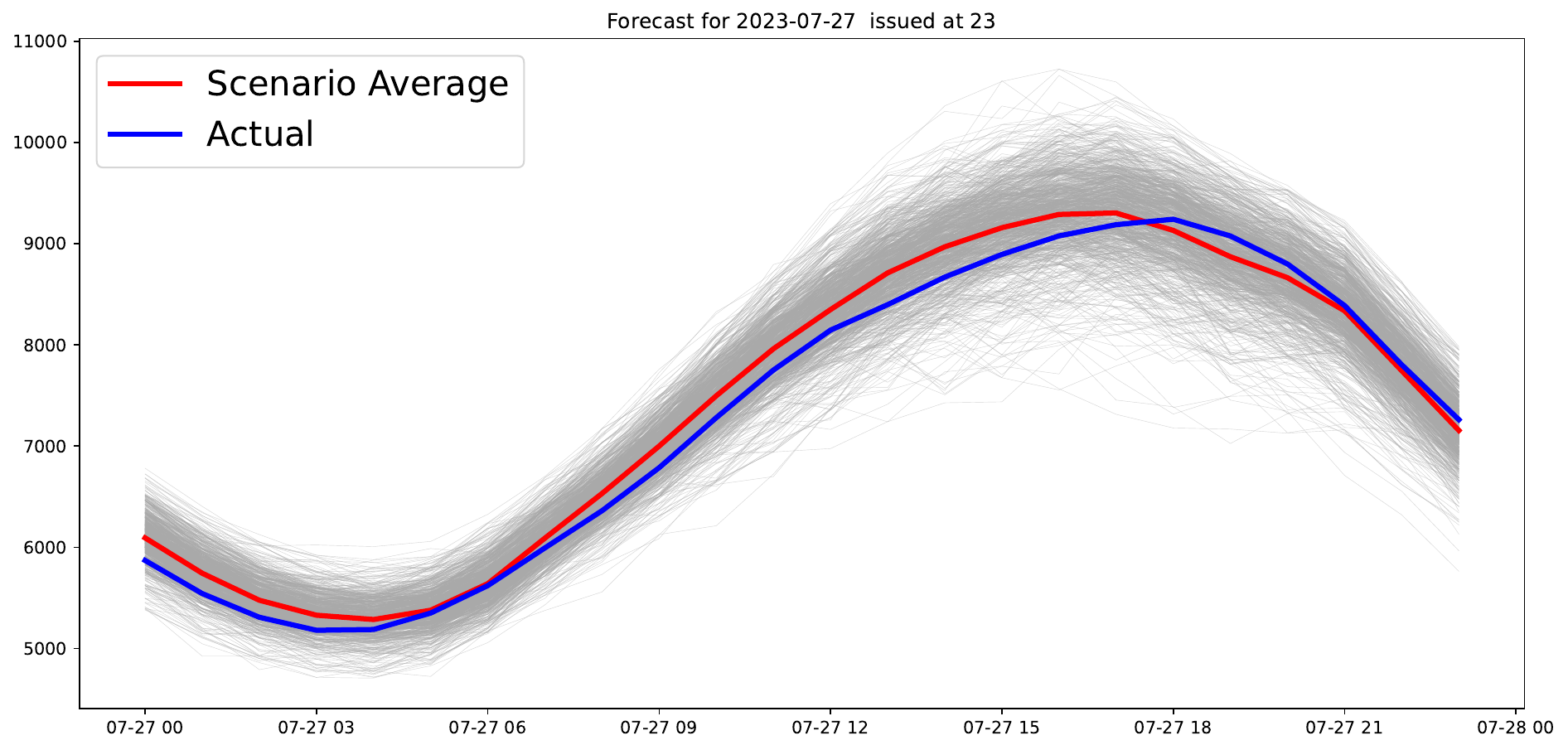}
\includegraphics[width=5cm]{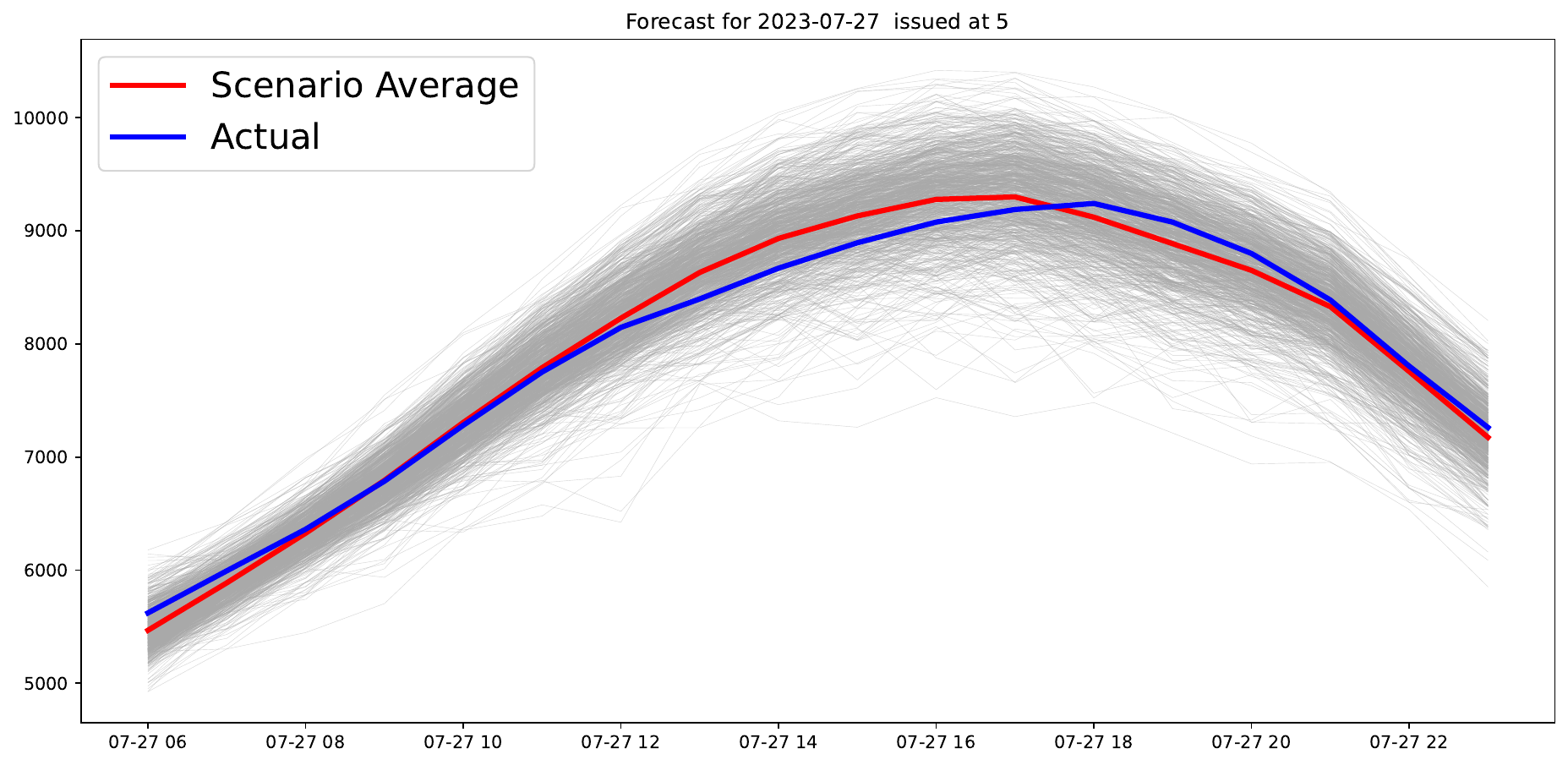}
\includegraphics[width=5cm]{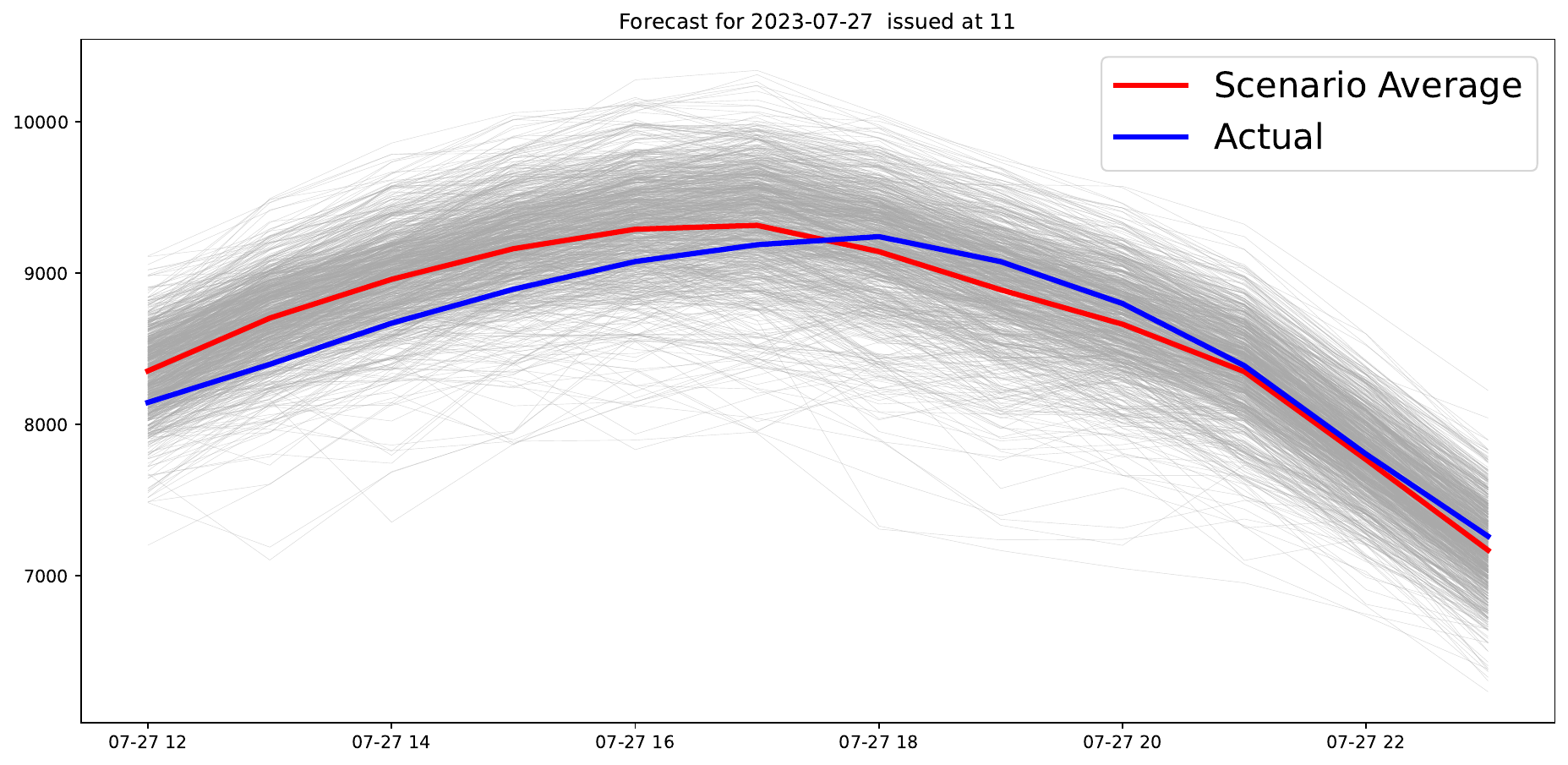}
\caption{PLProb scenarios for 2023-07-27 using forecasts generated at $23, 05, 11$ (in that order from left to right)}
\label{fig:scenario_pseg_1}
\end{figure}

\begin{figure}[htb!]
\includegraphics[width=5cm]{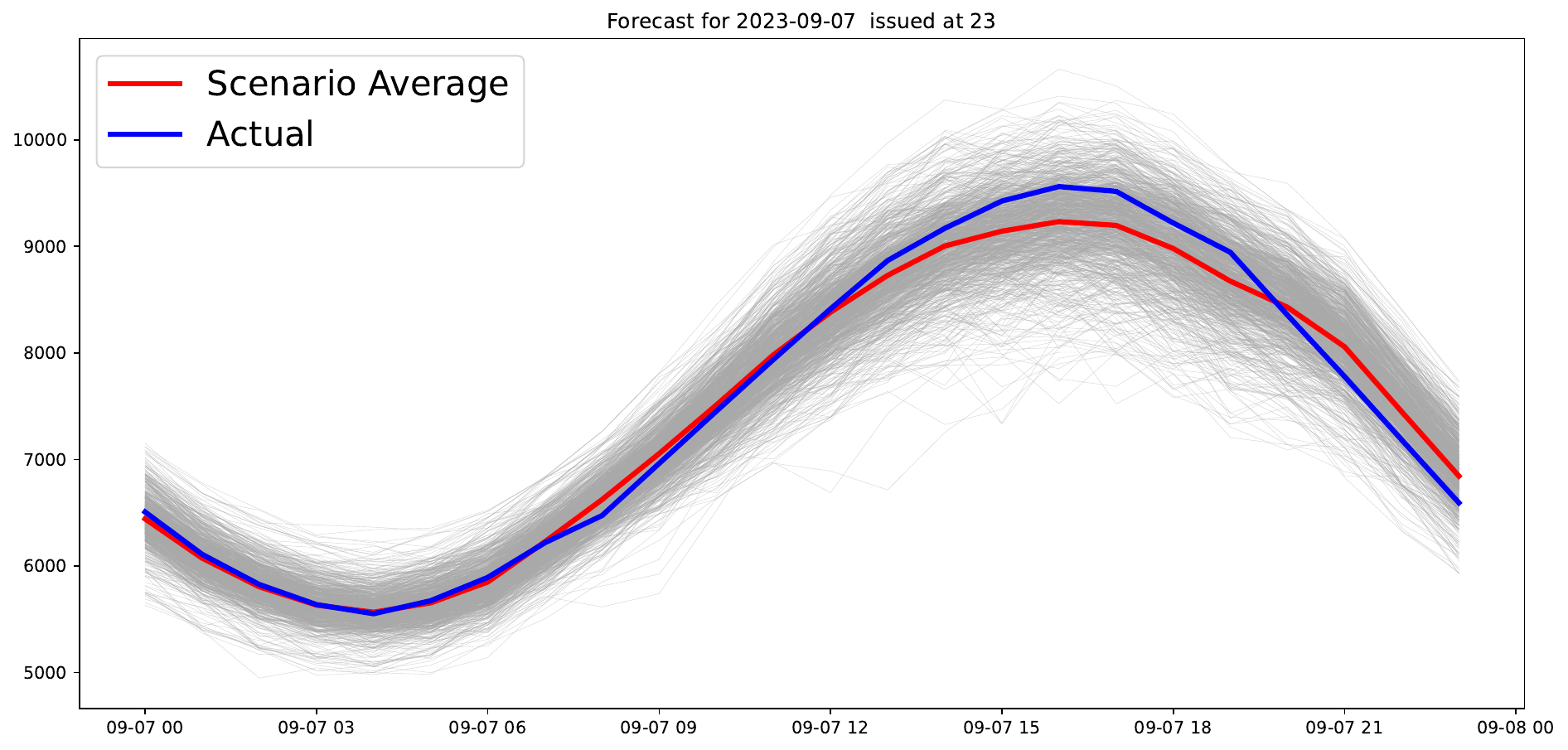}
\includegraphics[width=5cm]{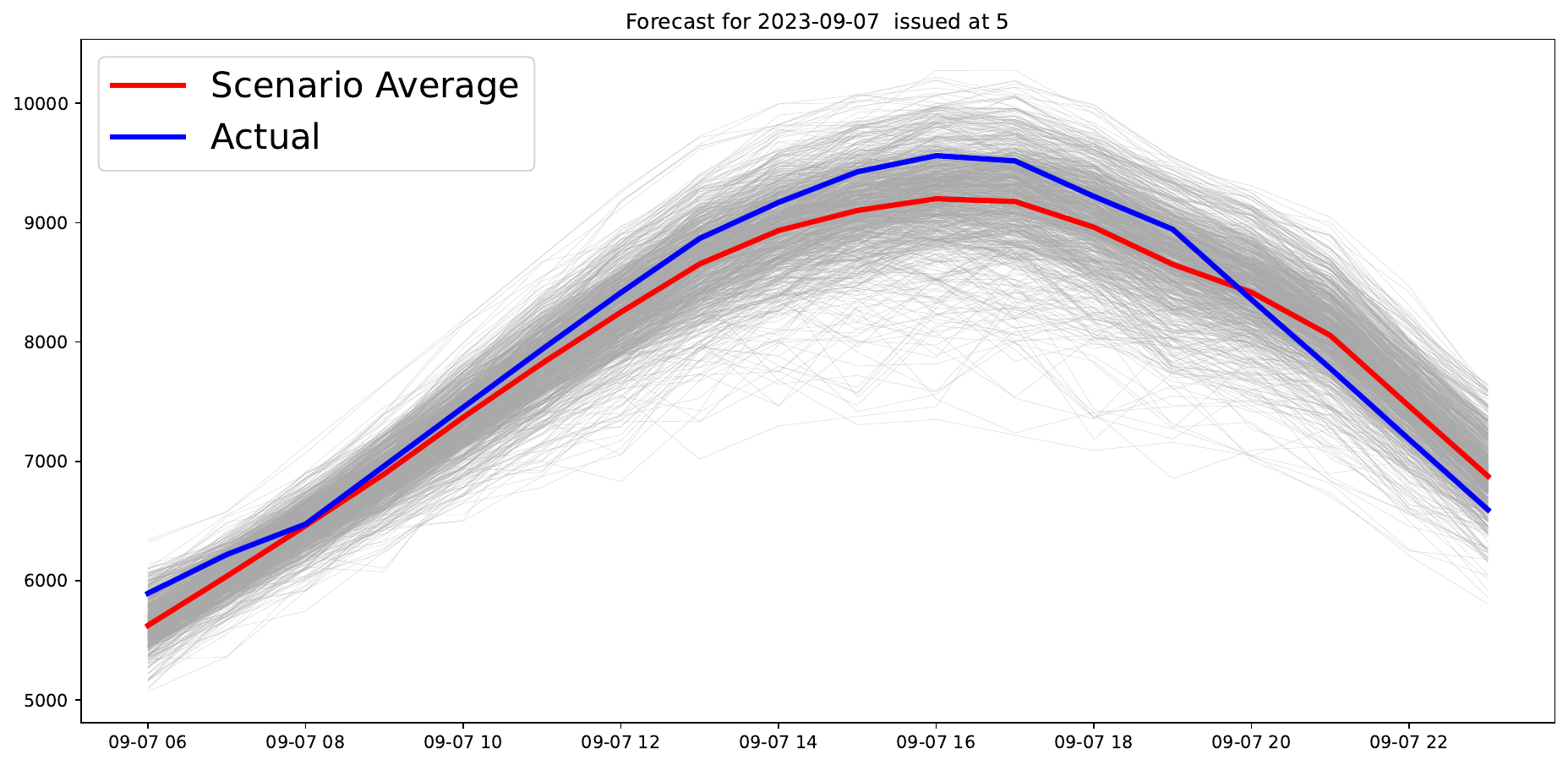}
\includegraphics[width=5cm]{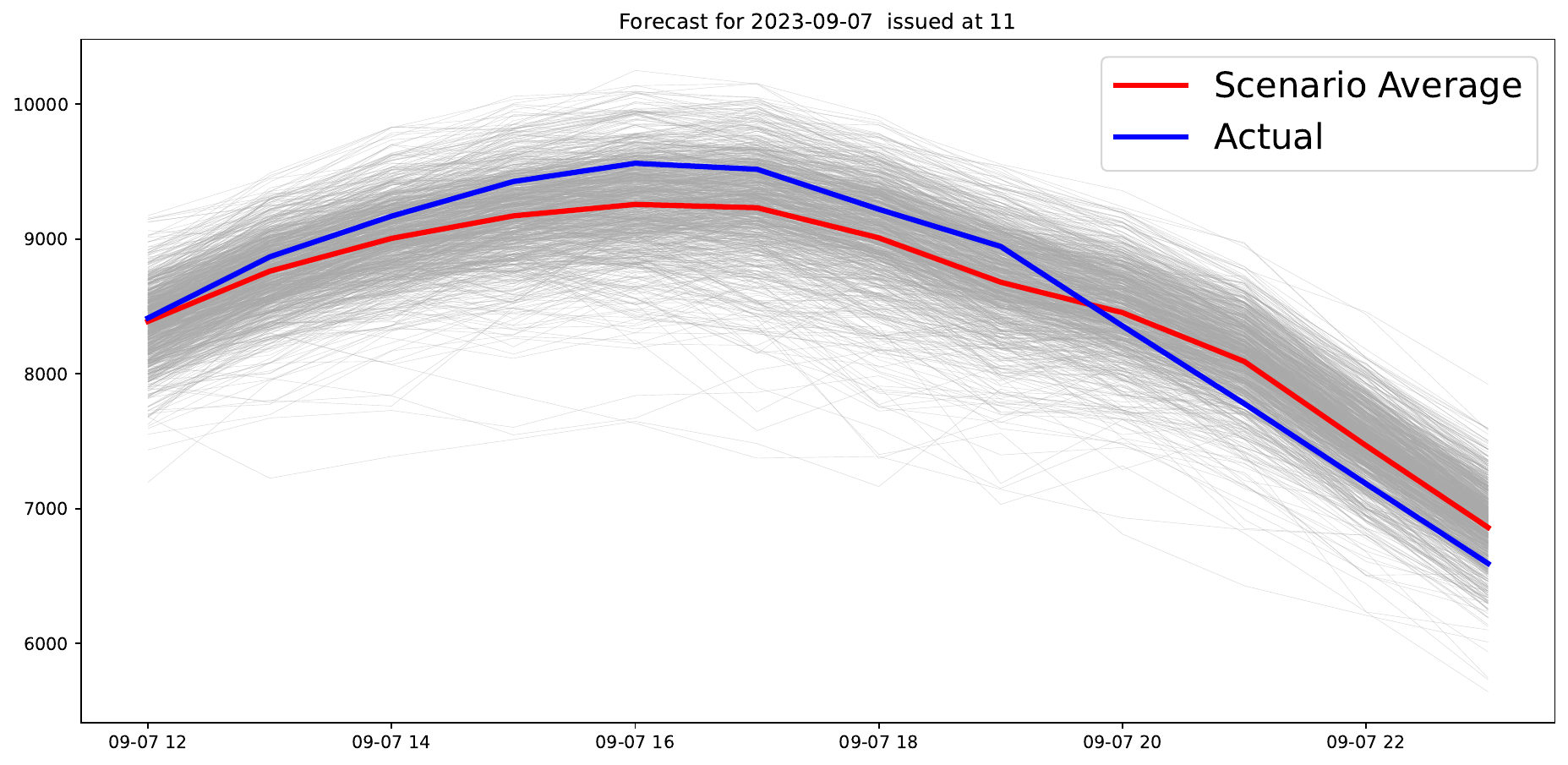}
\caption{PLProb scenarios for 2023-07-05 using forecasts generated at $23, 05, 11$  (in that order from left to right)}
\label{fig:scenario_pseg_2}

\end{figure}

Figures \ref{fig:scenario_pseg_1} and \ref{fig:scenario_pseg_2} present sample scenarios generated for different days given the forecast at three different times: 
\begin{itemize}
    \item the actual load (in blue) of PSE\&G;
    \item a batch of $K = 1000$ \texttt{PLProb} scenarios (in light grey) for PSE\&G;
    \item  the average of the scenarios (in red) computed with the conditional scheme based on the MIDATL forecasts labeled $23, 5, 11$. 
\end{itemize}

It is important to note that the current actual load is only plotted for reference, it is not used for the scenario generations. We only use past data for the fitting of the distributions described earlier. 

\subsection{Predicting Probability of Hitting the Running Maximum Load}   
We  look at the problem of estimating the probability that today is the new running CP day of the current year and adapt the estimators to our new setting. 

Let us denote by $d^*$ the current date, by $\{\bar{L}^d_{z,h}\}_{d \in \NN^*, h \in \{0, \dots, 23\}}$ the actual load of zone $z$ on day $d$ during the hour $h$ and by $\{M^d_{2}\}_{d \in \NN^*}$ the daily maximum actual load of zone 2 on day $d$. We denote by $\operatorname{CP}_{d^*} = \underset{d < d^*}{\operatorname{argmax}} \{ M^d_{2} \}$ the running maximum seen up until the day prior to $d^*$. Let $K$ be the number of generated scenarios and we denote by  $l^{d,k}_{z,h}$ the load of zone $z$ on day $d$ at hour $h$ and by $m^{d,k}_2$ the maximum daily load of scenario $k \in \{1, \dots, K\}$ for day $d$.

We compute the time-series of probabilities that today's peak is higher than the previous peaks by computing the frequency of scenarios whose daily maximum load exceeds the running maximum of days prior to $d^*$: 
$$ \widehat{\operatorname{prob}}_1 = \frac{1}{K} \sum_{k=1}^K \mathbbm{1}_{\{m^{d^*,k}_2 > \operatorname{CP}_{d^*}\}} $$
with the convention that $\operatorname{CP}_{0} = 0$. 

\begin{figure}[htb!]
    \begin{subfigure}[b]{\textwidth}
    \centering 
    \includegraphics[width=16cm]{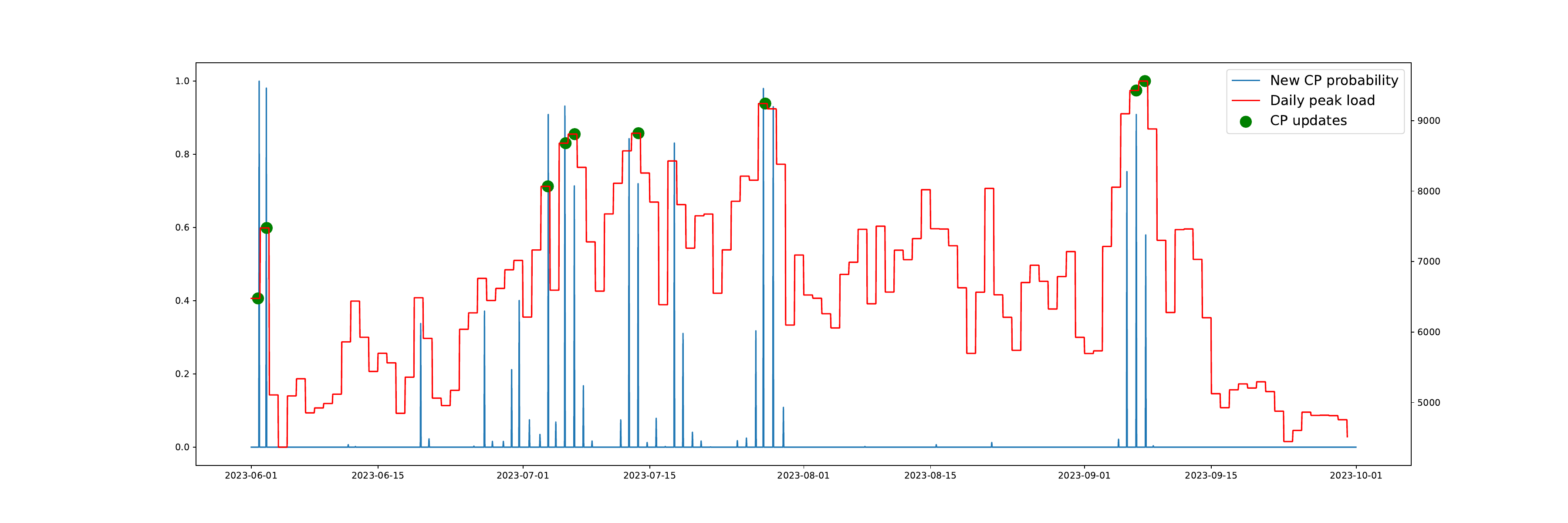}
    \caption{Using Forecast (23)}
    \end{subfigure}
    
\begin{subfigure}[b]{\textwidth}
\centering 
\includegraphics[width=16cm]{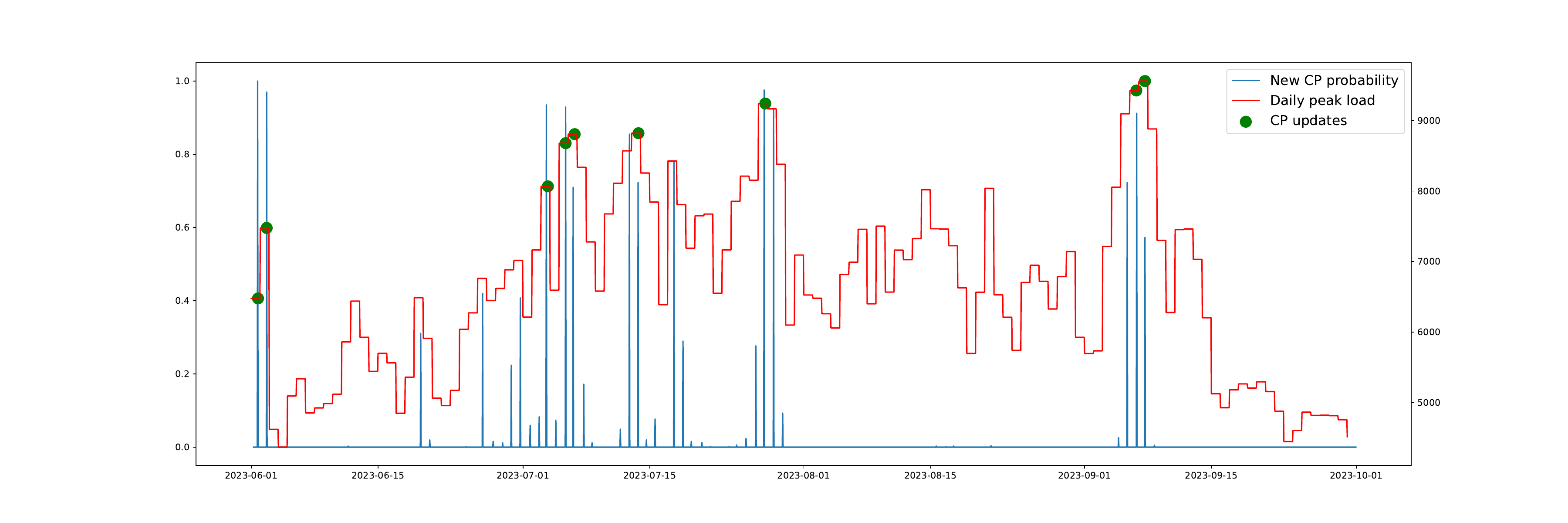}
        \caption{Using Forecast (05)}
\end{subfigure}
\begin{subfigure}[b]{\textwidth}
\centering 
\includegraphics[width=16cm]{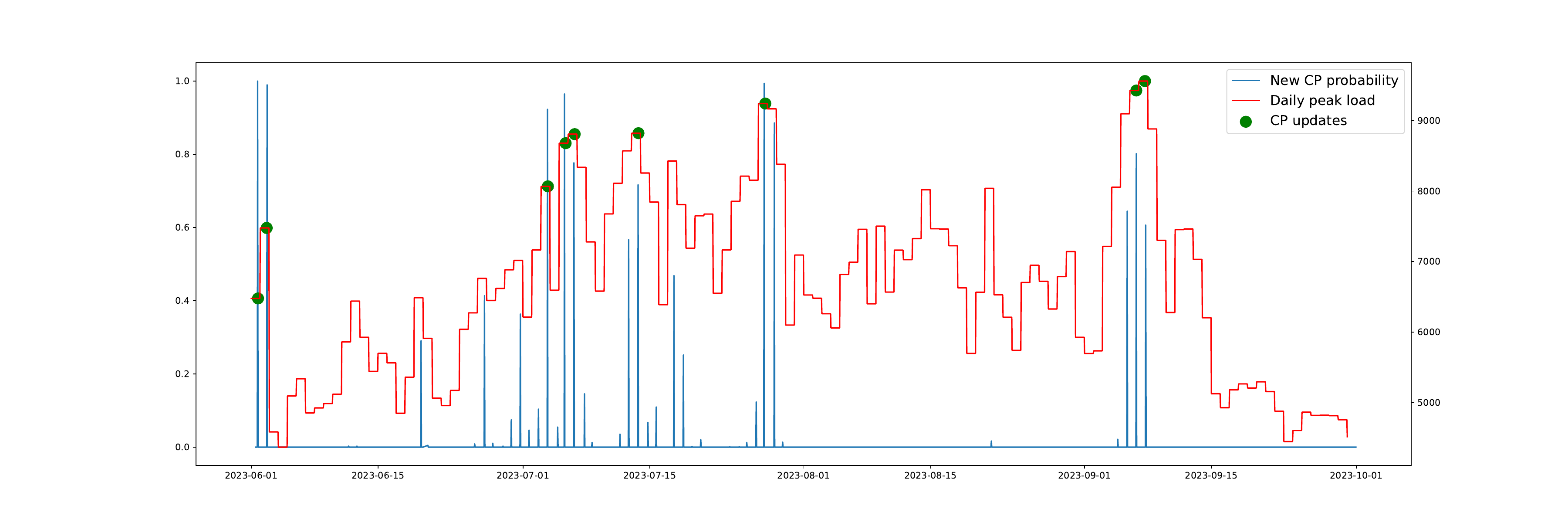}
        \caption{Using Forecast (11)}
\end{subfigure}
    \caption{Running new CP probability and maximum daily load time-series with different forecasts and a scenario count of 1000.}
    \label{fig:PSEG_proba_ts_frcst_2023}
\end{figure}

For each of the forecast sets $23, 05, 11$, we plot in Figure \ref{fig:PSEG_proba_ts_frcst_2023}:  
\begin{itemize}
    \item The actual daily PSE\&G peak load (red continuous curve)
    \item Green dots indicating the days when the daily maximum load exceeded the running maximum up until now;
    \item For each day, a blue vertical bar whose height is the probability that today's maximum load will exceed the previous ones.
\end{itemize}

Each blue spike corresponds to elevated probabilities that the days in question are coincident peak days. We see that this sequential method manages to predict high probabilities for all high daily maximum load days, but more importantly all the 8 running CP day updates (whenever a new running maximum is reached). The number of spikes is also fairly moderate (13 with a probability above 0.5). Section \ref{sec:backtest_implications} will analyze the performance of more involved strategies on a longer testing dataset. We also notice that using the most recent forecast leads to lower probabilities being assigned to False Positives, namely the days that do not see a CP update. This phenomenon is more noticeable on the graph using forecast (11) with lower spikes around July 01, 2023.

\subsection{Estimation of the time of the daily peak}

Given that a day is predicted to be a CP-day, we now need to estimate the probability of each hour being a CP-hour. We are looking for the probability that each hour $h$ sees the maximum daily load from the scenarios $k \in \{1, \dots, K\}$ on that day and we denote it by $\operatorname{prob}^{h}_1$. This leads to a vector of length $24$ containing the estimators: 
$$ \widehat{\operatorname{prob}}^{h}_1 = \frac{1}{K} \sum_{k=1}^K \mathbbm{1}_{\{l^{d^*, k}_{2,h} = m^{d^*,k}_2\}}, \quad h \in \{0, \dots, 23\}$$

We compute these estimates for the days on which the CP updates occur and plot: 
\begin{itemize}
\item the hourly actual load of the day as a continuous red curve,
\item the histogram of the hourly probabilities that a given hour is the hour the daily maximum occurs, computed from the generated scenarios,
\item a point in green indicating the hour when the maximum load actually occurs. 
\end{itemize}

\begin{figure}[htb!]
\begin{subfigure}[b]{0.45\textwidth}
\centering 
\includegraphics[width=\textwidth, height = 7cm]{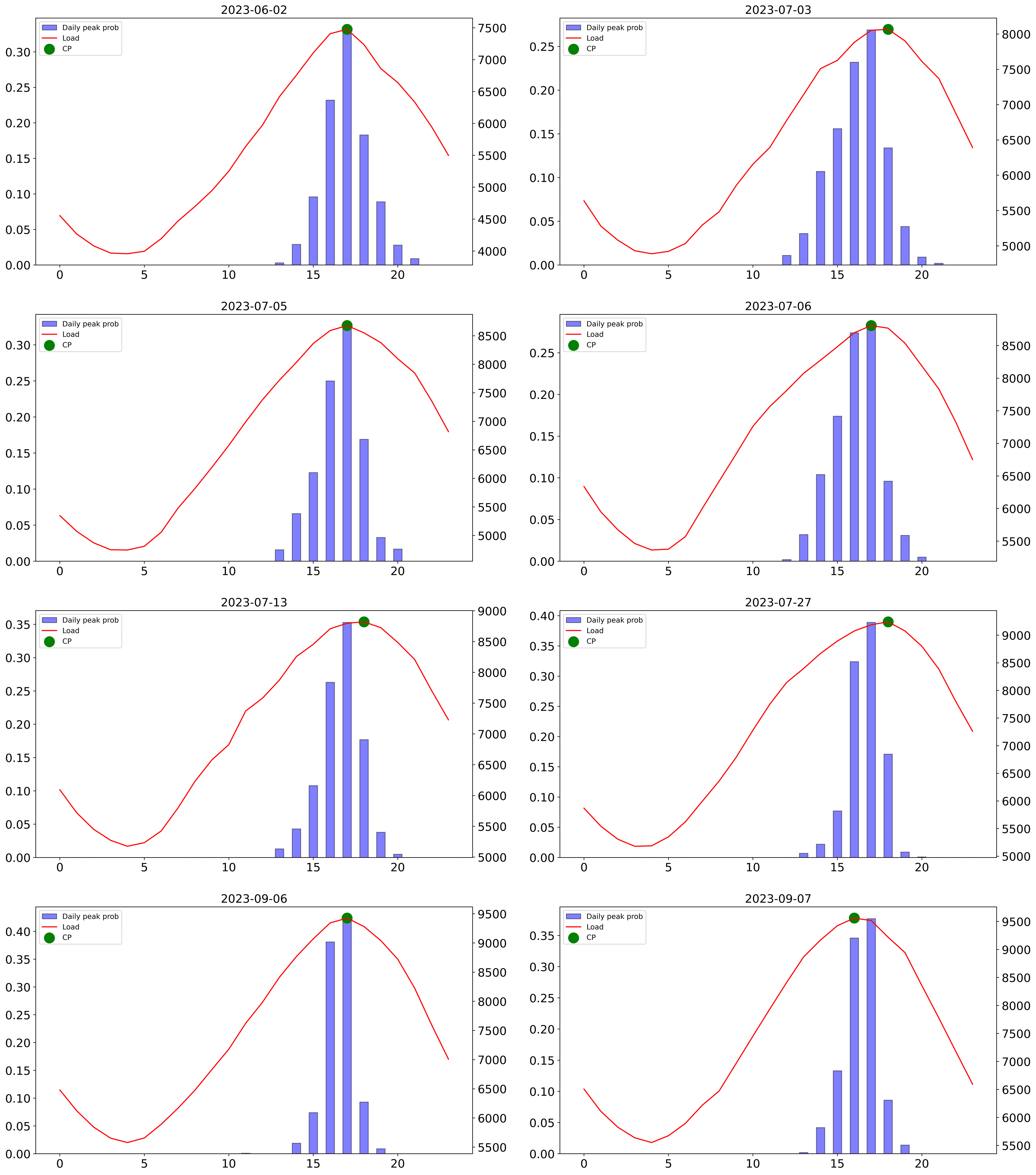}
\caption{Using forecast (23)}
\end{subfigure}
~ 
\begin{subfigure}[b]{0.45\textwidth}
\centering 
\includegraphics[width=\textwidth, height = 7cm]{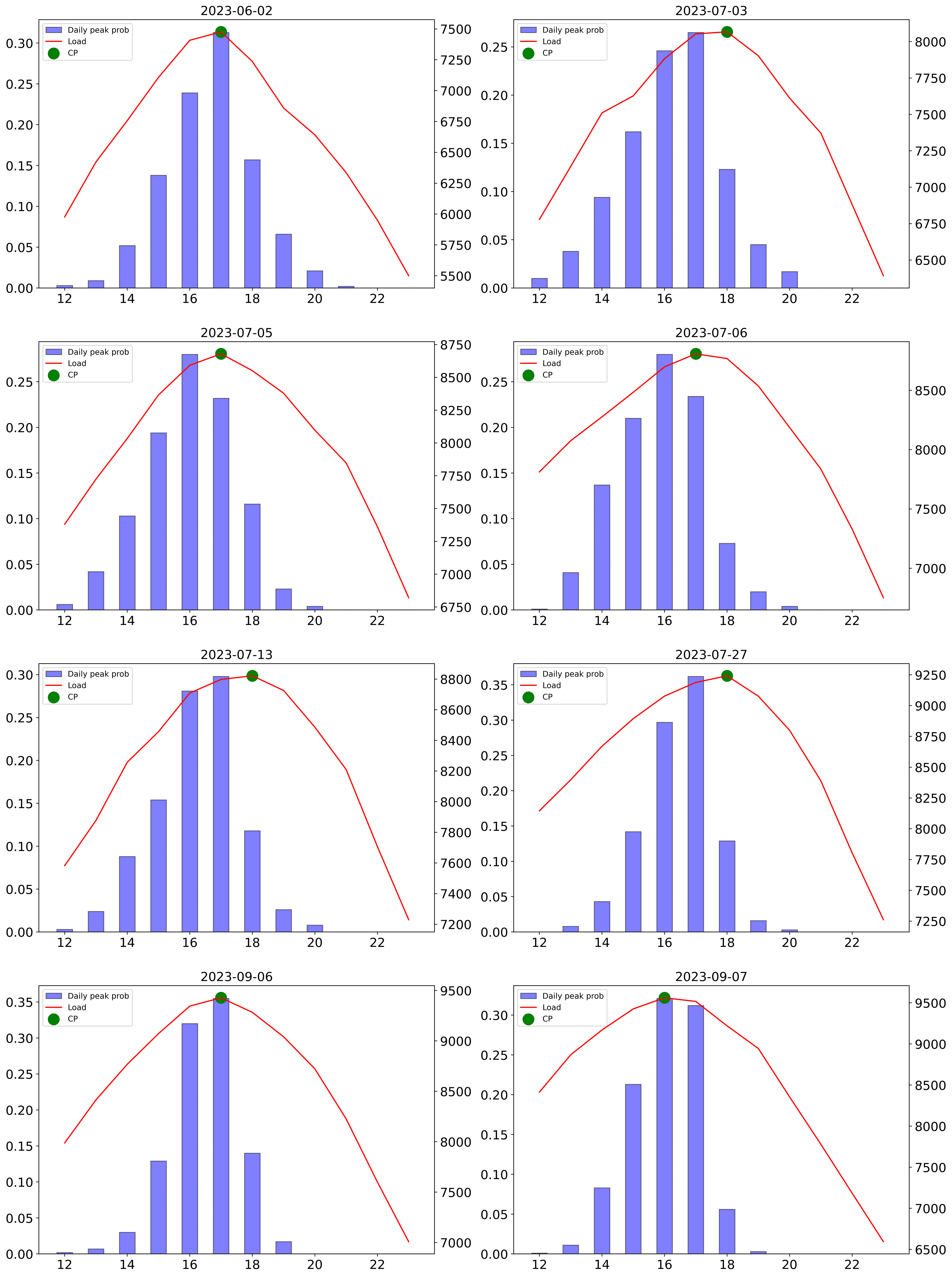}
\caption{Using forecast (11)}
\end{subfigure}
\caption{PSEG Daily Maximum hourly probability based on two forecasts.}
\end{figure}


\section{Extension to Multiple Coincident Peak Programs}
\label{sec:multi_cp}

This section is dedicated to the treatment of multiple CP programs in the same setting as in \ref{sec:frcst_iso} where both actual load and forecast datasets are available for the jurisdiction of interest. For the sake of illustration, we will focus on the 5CP program from PJM described in Subsection \ref{subsec:PJM_pseg_cp_program}.  We therefore use the corresponding notations. 

\subsection{Input of the Model}

We assume the availability of the following datasets: 
\begin{itemize}
    \item the actual load historical data of the RTO at an hourly resolution,
    \item the day-ahead forecast data of the RTO at an hourly resolution.
\end{itemize}
For instance, we consider the same load forecast data in Section \ref{sec:Conditional Simulation Results} with four different forecasting horizons. We let $h_s$ be the first hour of the day-ahead load forecast. For example, for the forecast published at 5:45 am, $h_s=6$.

\subsection{Method}

The problem of predicting 5 CPs is similar in nature to the 1CP case and requires to carry a list of the 5 prior CP events across time. We denote by $\text{CP}^*_i$ the running daily maximum load of rank $i$ reached at time $\tau_i$:

\begin{align*}
   & \text{CP}^*_1 = \underset{d < d^*}{\operatorname{max}} \{M_d\}, \quad \tau_1 = \underset{(d,h) : \, d < d^*, \, h \in \{h_s, \dots, 23\}}{\operatorname{argmax}} \{\overline{L}_{d,h}\}, \\ 
   & \text{CP}^*_k = \underset{d < d^*, d \notin \{\tau_1, \dots, \tau_{k-1}\}}{\operatorname{max}} \{M_d\}, \quad \tau_k = \underset{(d,h) : d < d^*, \, d \notin \{\tau_1, \dots, \tau_{k-1}, \, h \in \{h_s, \dots, 23\}\}}{\operatorname{argmax}} \{\overline{L}_{d,h}\},  \quad k \in \{2, \dots, 5\}.
\end{align*}

\begin{figure}[htb!]
\includegraphics[width=7cm]{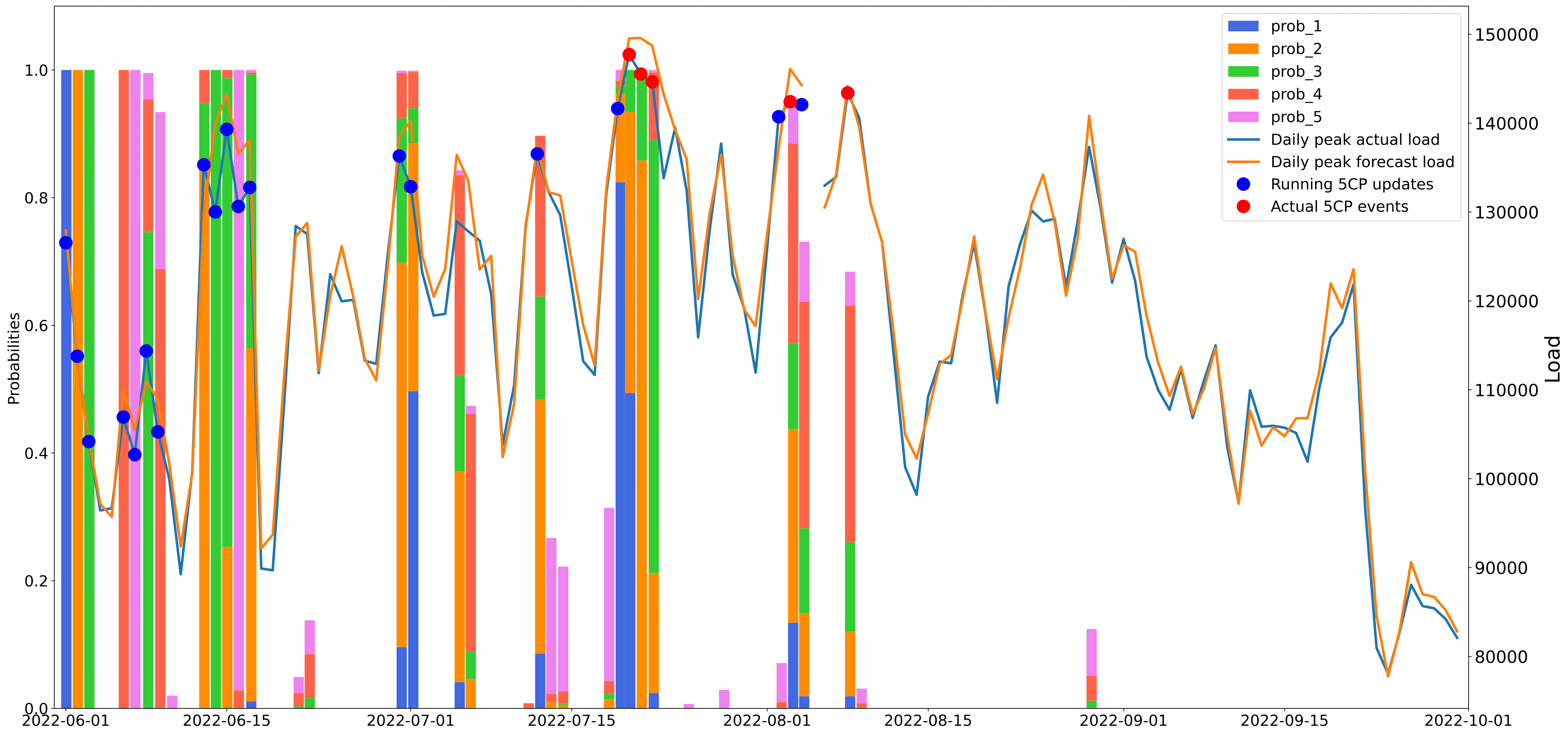}
\includegraphics[width=7cm]{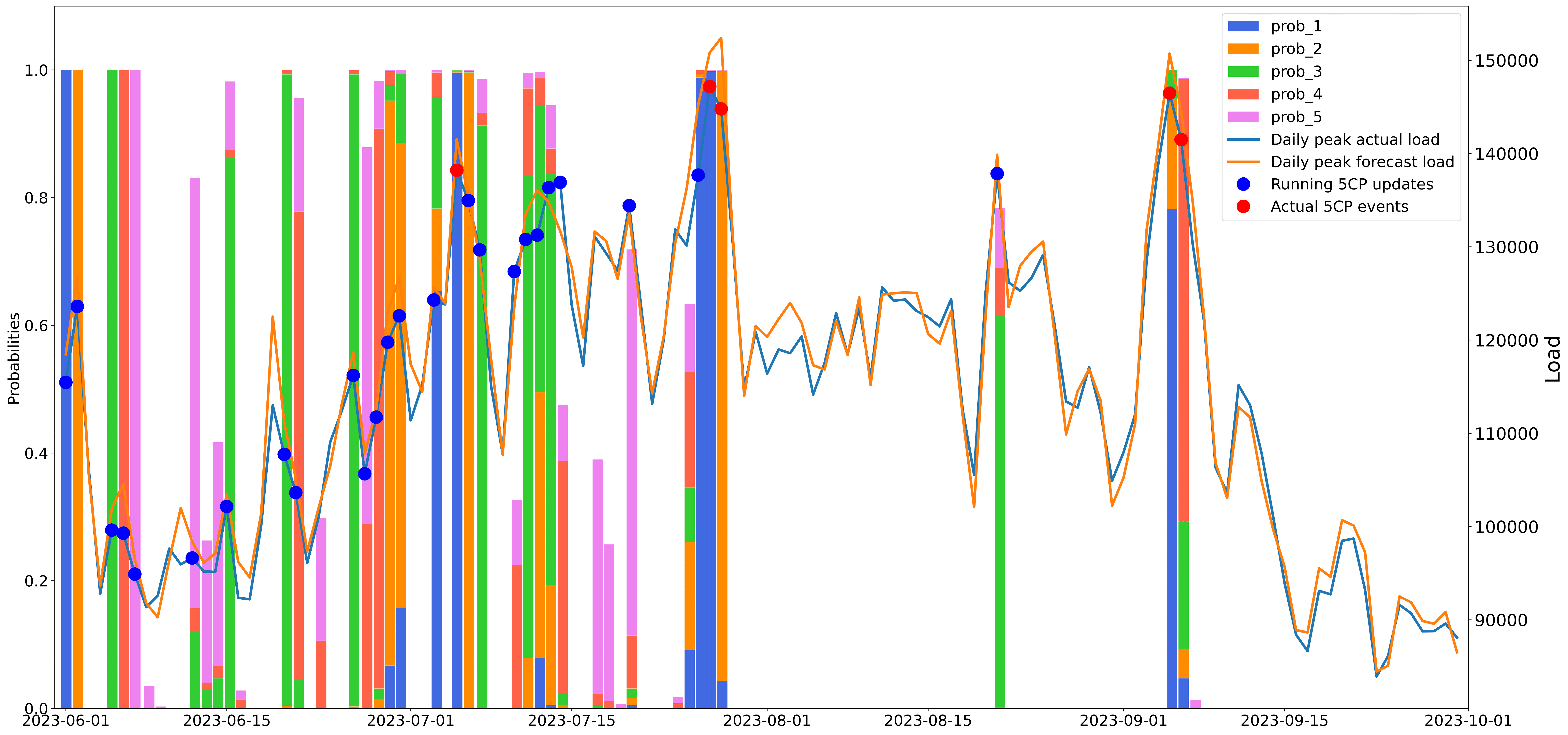}
\caption{Probabilities from forecasts (11) for the year 2022 (left pane) and 2023 (right pane)}
\label{fig:five_cp_2022-2023}
\end{figure}

We compute the time-series of five different probability estimates for :  
\begin{itemize}
    \item $\text{prob}_1 = \mathbb{P}(\text{daily maximum load} \geq \text{CP}^*_1)$,
    \item $\text{prob}_k =\mathbb{P}(\text{CP}^*_{k-1} > \text{daily maximum load}  \geq \text{CP}^*_k), \quad k = 2, \dots, 5$.
\end{itemize}
where $\text{CP}^*_i, i \in \{1, \dots, 5\}$ is the running CP value of rank $i$ up until $d^*$ excluded. We are interested in the probability that a day is a potential CP day and hence compute the sum of these probabilities, which are represented  by the stack bars in the graphs from Figure \ref{fig:five_cp_2022-2023}. 

\begin{remark}
The missing point in the year 2022 corresponds to missing data for 2022-08-05 in the database of PJM, it does not affect the final CP values. Similarly as in the 1 CP case, the first five spikes are elevated because of the need to populate the list of our running 5 CP events. 
\end{remark}

The main issue with this approach is that there are too many warnings, including the first five peaks (while they are most likely not going to be any of the five CPs). The goal is then to adopt a data-based approach to trim down most of the False Positives (elevated probabilities for days that do not end up belonging to any of the true 5 CPs). 
    
\subsection{Threshold 5CP event probabilities}
\label{subsec:thresh_5cp} 

We compute a data-driven percentile threshold \texttt{thresh} which corresponds to the $k$-th percentile of the past daily maximum load for $k \in \{ 80, 90, 95\}$. This selection is motivated by the fact that there are about 80-88 business days in a PJM summer period and we are interested in the 5 top values, which represent around 6\% of the days. 95 might be too stringent and we may miss some values, while 80 might be too large and give too many warnings. 
We define modified running coincident peaks, denoted as $\widetilde{\text{CP}}^*_k = \max(\text{CP}^*_k, \texttt{thresh})$. 
We compute the following 5 time-series of different probabilities: 
\begin{itemize}
    \item $\text{prob}_1 = \mathbb{P}(\text{daily maximum load} \geq \widetilde{\text{CP}}^*_1)$,
    \item $\text{prob}_k =\mathbb{P}(\widetilde{\text{CP}}^*_{k-1} > \text{daily maximum load} \geq \widetilde{\text{CP}}^*_k), \quad k = 2, \dots, 5$,
\end{itemize}
where $\text{CP}_i, i \in \{1, \dots, 5\}$ is the running CP value of rank $i$. Table \ref{tab:thresh_values} provides the cutting values for the chosen percentiles. 

\begin{table}[htb!]
    \centering
    \begin{tabular}{|c|c|c|}
    \hline
        Percentile & Threshold used for 2022 &  Threshold used for 2023\\
        \hline 
        95 & 137,831&137,405\\
        90 &132,349 & 132,048\\
        80 & 125,182& 125,098\\
    \hline 
    \end{tabular}
    \caption{Past data $k$-th percentile}
    \label{tab:thresh_values}
\end{table}

\begin{figure}[htb!]

\includegraphics[width=7cm]{5CP_figures/cp5_thresh_figs/fig1_fill_2022_thresh_0.png}
\includegraphics[width=7cm]{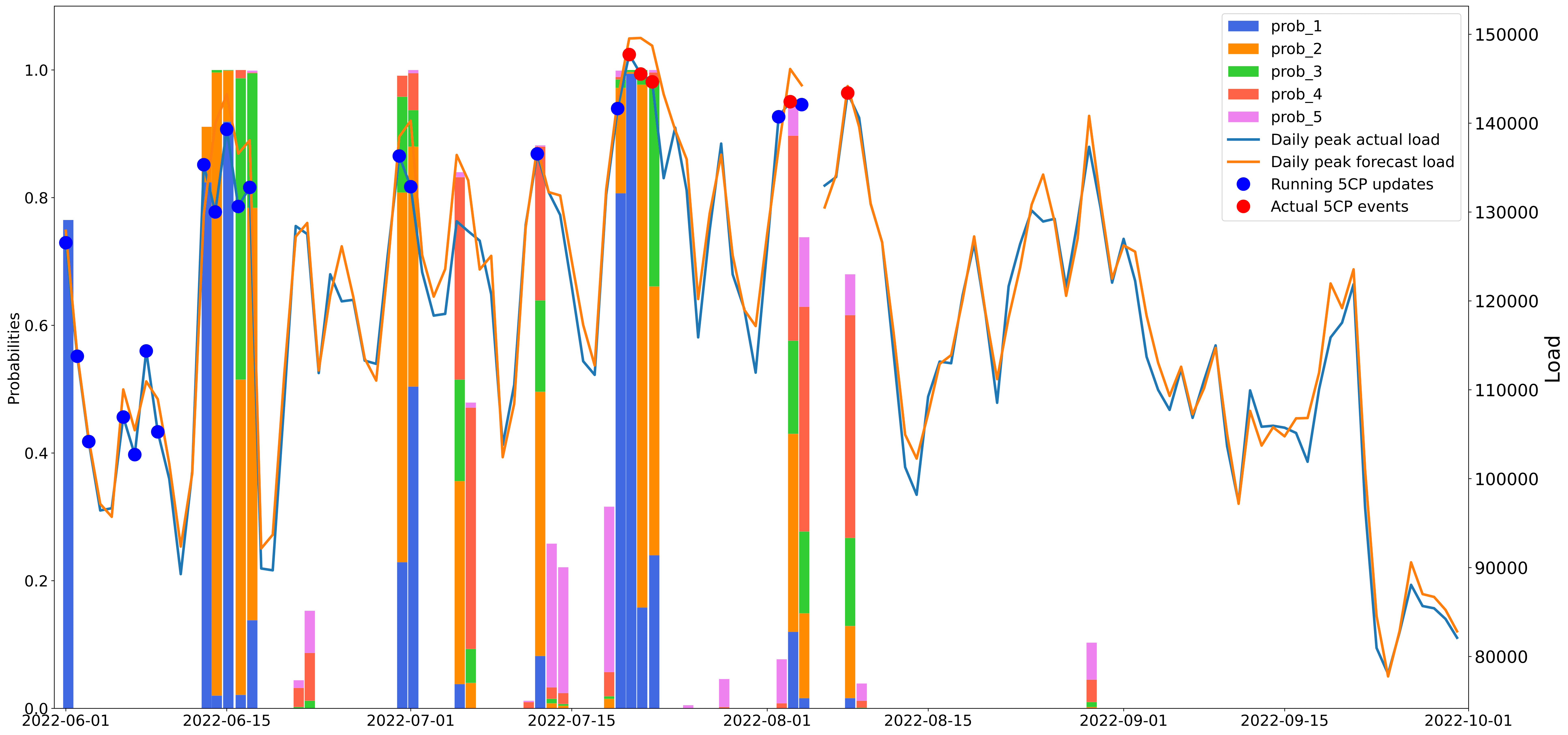}
\includegraphics[width=7cm]{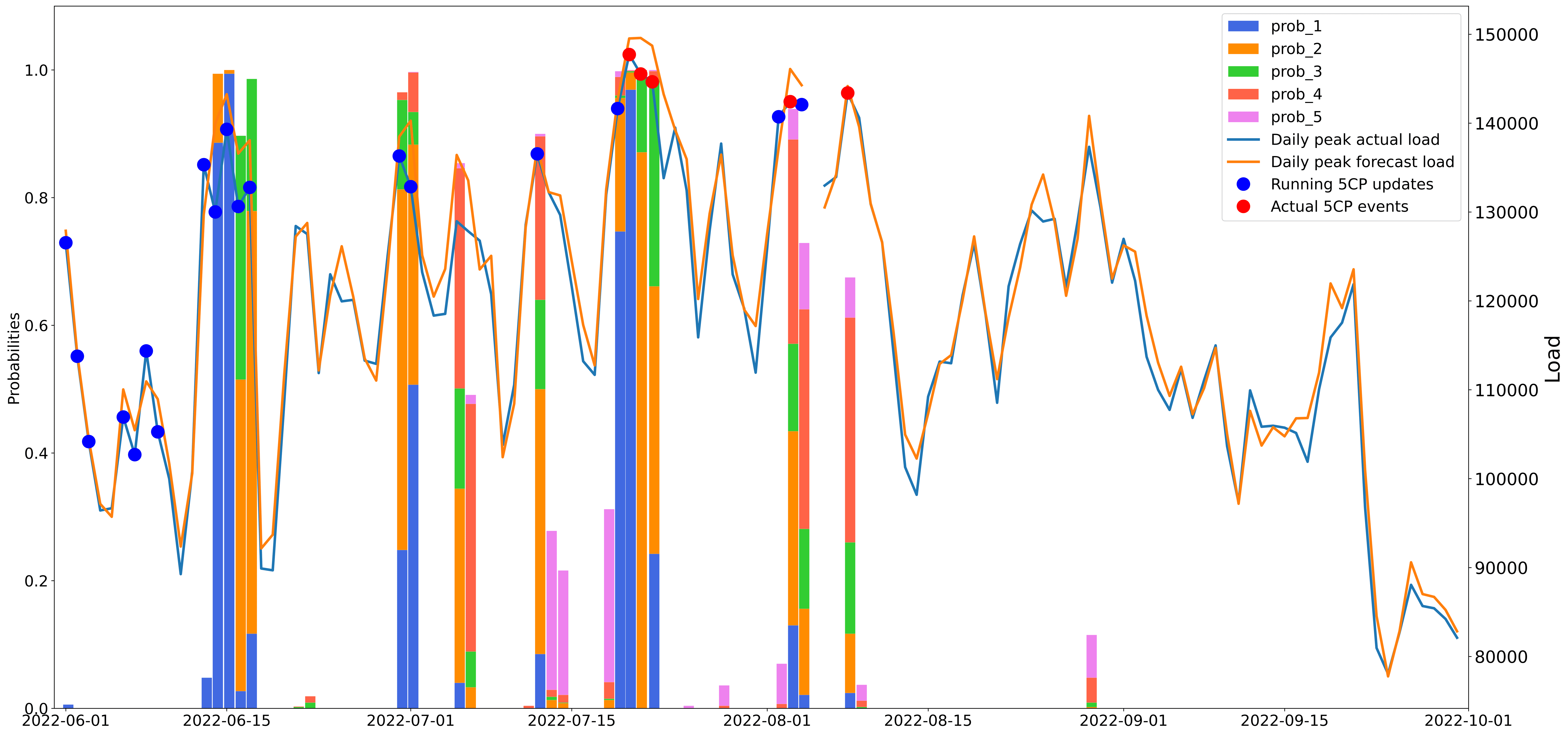}
\includegraphics[width=7cm]{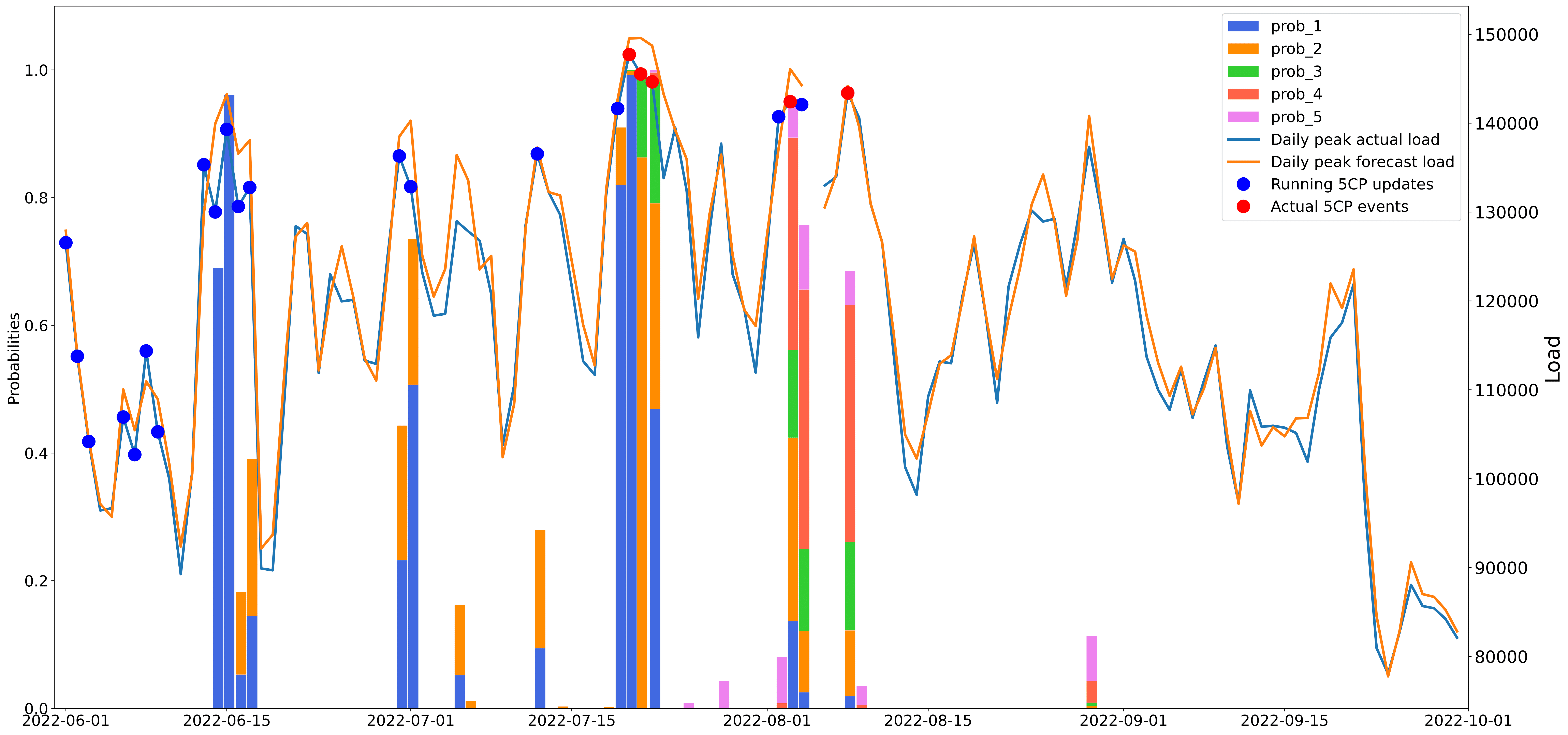}

\caption{Probabilities from forecasts at 11:00 for 2022 with thresholds of $k$-th percentile for $k = \text{None}, 80, 90, 95$ (respectively from top-left to bottom-right)}
\end{figure}

\begin{figure}[htb!]

\includegraphics[width=7cm]{5CP_figures/cp5_thresh_figs/fig1_fill_2023_thresh_0.png}
\includegraphics[width=7cm]{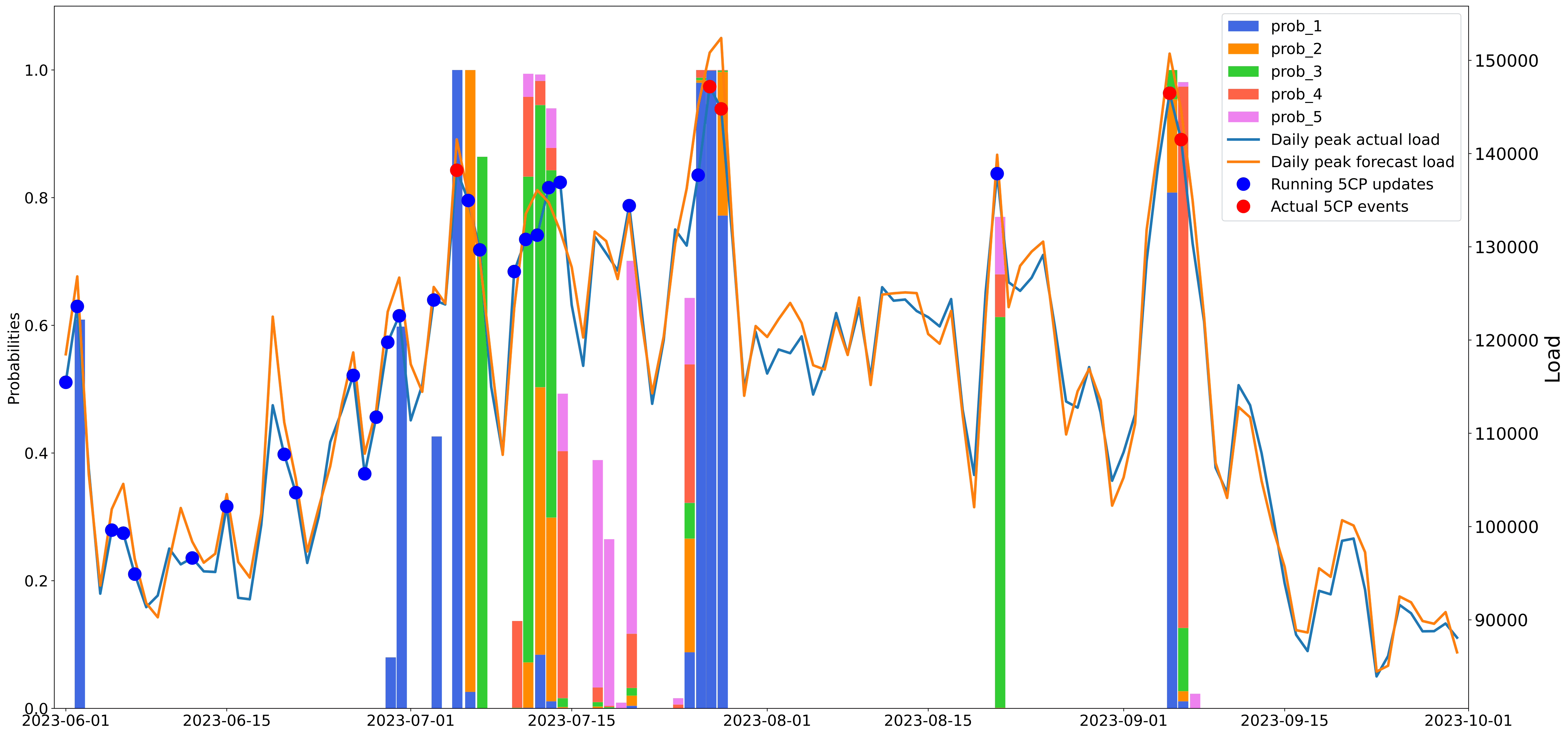}
\includegraphics[width=7cm]{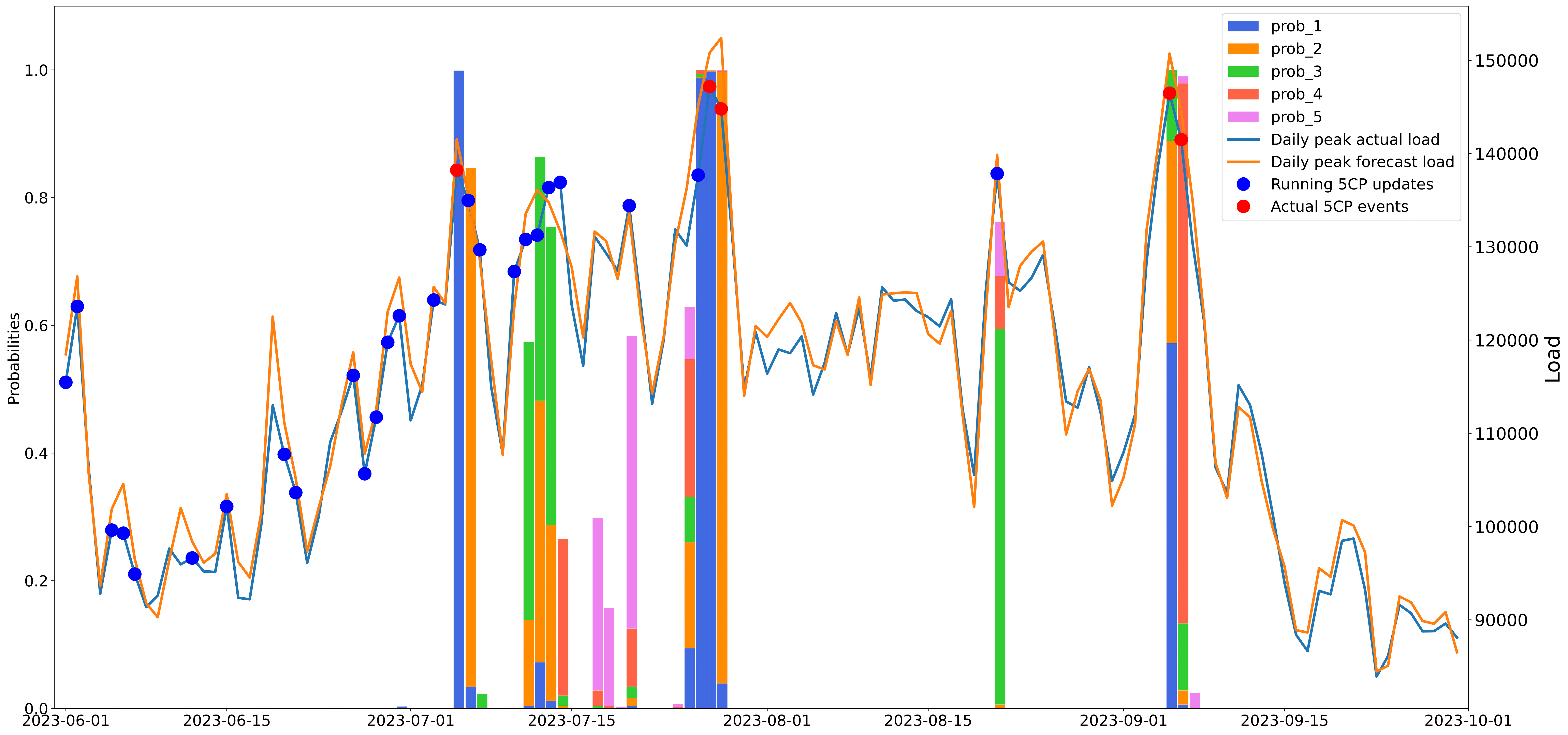}
\includegraphics[width=7cm]{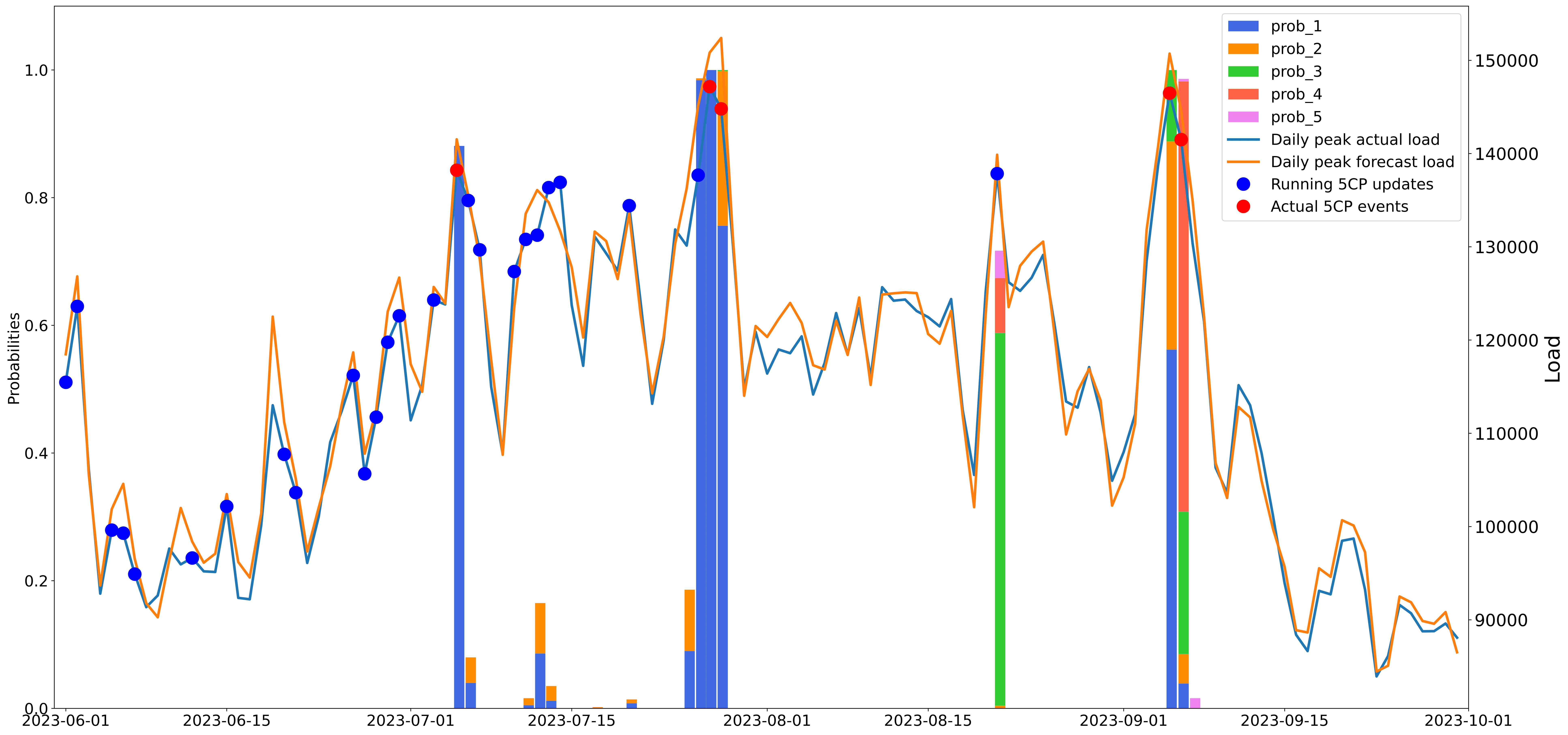}

\caption{Probabilities from forecasts at 11:00 for the year 2022 with thresholds of $k$-th percentile for $k = \text{None}, 80, 90, 95$ (respectively from top-left to bottom-right)}
\end{figure}

\begin{table}[htb!]
    \centering
    \begin{tabular}{|c|c|c|}
    \hline
        Percentile & $\# \{ \sum_{i=1}^5 prob_i > 0.5\}$ (5CP) - 2022 &  $\# \{ \sum_{i=1}^5 prob_i > 0.5\}$ (5CP) - 2023\\
        \hline 
        95 & 10 (5) & 7 (5)\\
        90 & 15 (5) & 13 (5)\\
        80 & 17 (5) & 16 (5)\\
        None & 29 (5) & 23 (5) \\
    \hline 
    \end{tabular}
    \caption{Warnings based on probabilities computed with threshold}
    \label{tab:warning_stats}
\end{table}

As expected, the lower the threshold the more warnings are sent. The choice of $k=90$ gives approximately 15 elevated probabilities of a new CP event. However, it is interesting to see that even with a more stringent threshold $k = 95$, we do not miss any of the final 5 CP events.

\section{Backtest and Load Curtailment Considerations} \label{sec:backtest_implications}

Motivated by the qualitative analysis of Subsection \ref{subsec:thresh_5cp}, we formulate different strategies and analyze their performance over ten years worth of data. We also draw some conclusions regarding the implementation of load curtailment, including the sizing or the deployment duration of the technological solutions. 

\subsection{Methodology}

We backtest the performance of our predictors over the period 2014-2023 (10 years in total) for the each of the considered CP programs. To that end, we split the corresponding dataset according to a sliding-expanding window scheme: given a year of interest $y^* \geq 2014$,  we divide the dataset in the following:
\begin{itemize}[widest={{Training set}.},leftmargin=*]
    \item[\textit{Training set}:] Business days in the summer period of the years between 2011 and $y^*$ (excluded), 
    \item[\textit{Testing set}:] Business days in the summer period of the year $y^*$. 
\end{itemize}

We define several signal generation strategies that we label by the threshold: 
\begin{itemize}[widest={{Strategy k}.},leftmargin=*]
    \item[\textit{Threshold} \texttt{1}.] The threshold is $0$. This corresponds to the naive method where we progressively populate the running coincident peak lists with only the values seen in the current year. 
    \item[\textit{Threshold} \texttt{2}.] The threshold is the past data $95^{\text{th}}$ percentile. 
    \item[\textit{Threshold} \texttt{3}.] The threshold is the past data $90^{\text{th}}$ percentile. 
\end{itemize}

We refine the above strategies by defining modified running coincident peak levels, denoted as $\widetilde{\text{CP}}_{d^*} = \max( \alpha \cdot \text{CP}_{d^*}, \text{thresh})$ and consider different versions depending on the value of $\alpha$: 
\begin{itemize}[widest={{Version k}.},leftmargin=*]
        \item[\textit{Version} \texttt{a}.] $\alpha = 1.0$ or 100\% CP lower 
        \item[\textit{Version} \texttt{b}.] $\alpha = 0.975$ or  97.5\% CP lower
        \item[\textit{Version} \texttt{c}.] $\alpha = 0.95$ or  95\% CP lower
        \item[\textit{Version} \texttt{d}.] $\alpha = 0.90$ or 90\% CP lower
    \end{itemize}

Hence, the method will lead to computing the survival function of the next day maximum load (i.e. the complement to 1 of its cumulative distribution function) at the point $\widetilde{\text{CP}}_{d^*}$: the daily maximum load will be compared to this threshold that is progressively adjusted as more CP updates are seen. 
In particular, versions  \texttt{a},  \texttt{b} and  \texttt{c} will be applied to the ISO CP programs, i.e. NYISO and PJM, while versions  \texttt{a},  \texttt{c} and  \texttt{d} to the utility CP program, i.e. PSEG. 

\vspace{\baselineskip}

Finally, we consider two types of signals for CP day classification: a naive one labeled ``Simple'' relying on the rule $ \widehat{\text{prob}}_1 \geq  0.5 $ and a categorical one labeled ``Color'' and based on the value of the probability: 
\begin{itemize}[widest={{Color Orange}.},leftmargin=*]
        \item[\textit{Simple} \texttt{S}.] Probability above 0.5,
         \item[\textit{Color} \texttt{C}.] Probability above 0.2 for NYISO and PJM, and probability above 0.4 for PSEG with more specific alerts broken down into:
        \begin{itemize}[widest={{Orange O}.},leftmargin=*]
        \item[\textit{Red} \texttt{R}.] Probability above 0.8,
        \item[\textit{Orange} \texttt{O}.]  Probability between 0.6 and 0.8,
        \item[\textit{Yellow} \texttt{Y}.]  Probability between 0.4 and 0.6.
        \item[\textit{Green} \texttt{G}.]  Probability between 0.2 and 0.4.
        \end{itemize}
    \end{itemize}
The idea behind the color signal \texttt{C} is to assign a priority order depending on the strength of the probability: hence, a red alert would correspond to days that have a very high probability of seeing a new CP,  an orange alert for a relatively high probability, while the yellow and green ones would be for predictions that are more in the gray area with a load that is expected to be quite close to the current maximum. 

We label the strategies as the concatenation of a threshold, a version and a signal type. Hence, strategy \texttt{1aS} corresponds to using threshold \texttt{1}, version \textit{a} and signal \texttt{S}:  this will constitute our benchmark. Other variants are possible by considering more conservative and/or sophisticated filtering rules. In the following subsections, we analyze the performance of the chosen strategies and derive some implications for practical load curtailment for capacity and transmission charge reduction.  

\subsection{Performance on Coincident Peak Day Prediction}

We summarize the performance of predicting coincident peak days for each ISO case study. The average performance across the testing period is summarized in Tables \ref{tab:nyiso_perf_cp_day}, \ref{tab:all_pseg_perf_cp_day_stats}, and \ref{tab:pjm_perf_cp_day} for the NYISO 1CP program, PSEG 1CP program, and PJM 5CP program respectively. The total number of alerts broken down per year is provided for a selection of combinations of strategy, version and signal type. 


\subsubsection{NYISO}
From Table \ref{tab:nyiso_perf_cp_day}, we first note that all strategies are performing well  in the prediction of the true single coincident peak per year, with a minimum average over 10 years of at least 0.9 Hence, the strategies miss at most only 1 CP with the simple signal. The under-performing strategies are the versions \texttt{a} with no thresholding, while versions \texttt{b} and \texttt{c} detect all of the 10 true CPs. This can be explained by the fact that the estimator will attribute lower probabilities to load that are forecasted to exceed the current CP by only a marginal amount, while in fact any load forecasted to be within the same range as the current running CP should be considered. 

\begin{table}[htb!]
    \centering
    \begin{tabular}{||c||c|c||c|c||}
    \hline 
       \multirow{2}{*}{Strategy}  & \multicolumn{2}{c||}{Simple (\texttt{S})}  &   \multicolumn{2}{c||}{Color (\texttt{C})} \\
       \cline{2-5} 
         & \# alerts &  \# CP &  \# alerts (\texttt{R}/\texttt{O}/\texttt{Y}) & \# CP (\texttt{R}/\texttt{O}/\texttt{Y})\\
       \hline 
        \texttt{1a} & 4.9 & 0.9 & 5.3 (3.5/0.9/0.9) & 1.0 (0.4/0.4/0.2) \\
        \texttt{1b} & 7.4 & 1.0 & 8.7 (5.0/2.0/1.7) & 1.0 (1.0/0.0/0.0) \\ 
        \texttt{1c} & 10.9 & 1.0 &  11.6 (8.1/2.4/1.1) & 1.0 (1.0/0.0/0.0)  \\
        \texttt{2a} & 4.3 & 0.9 & 4.6 (2.6/1.2/0.8) & 1.0 (0.4/0.4/0.2) \\
        \texttt{2b} & 6.6 & 1.0 & 7.8 (4.1/1.9/1.8) & 1.0 (1.0/0.0/0.0) \\ 
        \texttt{2c} & 9.9 & 1.0 & 10.6 (7.0/2.4/1.2) & 1.0 (1.0/0.0/0.0) \\
        \texttt{3a} & 4.3 & 0.9 & 4.6 (2.6/1.1/0.9) & 1.0 (0.4/0.4/0.2) \\
        \texttt{3b} & 6.6 & 1.0 & 7.7 (4.1/1.9/1.7) & 1.0 (1.0/0.0/0.0)\\
        \texttt{3c} & 9.9 & 1.0 & 10.9 (7.1/2.4/1.4) & 1.0 (1.0/0.0/0.0) \\
    \hline 
    \end{tabular}
    \caption{Average performance statistics of selected strategies for the 1CP program of NYISO. The averages are computed over 10 years, the color strategy shows the overall average performance broken down according to the three colors: red \texttt{R}, orange \texttt{O}, and yellow \texttt{Y})}
    \label{tab:nyiso_perf_cp_day}
\end{table}

\begin{figure}[htb!]
    \centering
    \includegraphics[width=16cm]{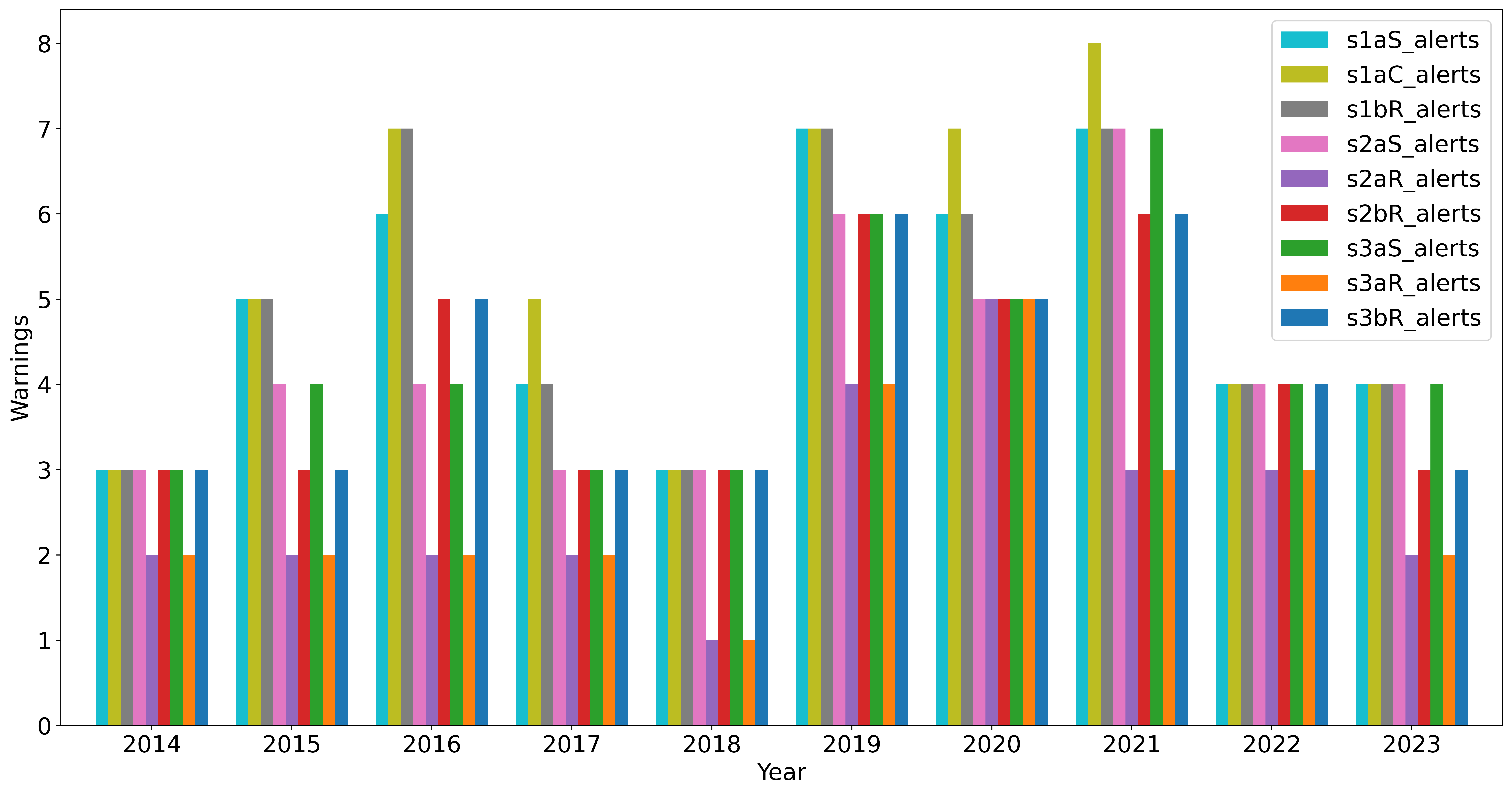}
    \caption{Breakdown per year of the number of alerts generated by selected strategies on NYISO 1CP program}
    \label{fig:nyiso_perf_cp_day_year}
\end{figure}

From the yearly breakdown in Figure \ref{fig:nyiso_perf_cp_day_year}, we see that there are more alerts associated with lower thresholding method: for instance, as many as 8 for year 2021 for the color signal \texttt{C}), while strategies \texttt{2aR} and \texttt{3aR} have the least number of warnings but only catch a small fraction of the true CPs, totaling 4 out of the 10 CPs. A trade-off has to be found in order to both maximize the number of caught CPs and minimize the number of false alerts.


\subsubsection{PSEG}

From Table \ref{tab:all_pseg_perf_cp_day_stats}, we first note that all strategies using no threshold (threshold = 0) lead to missing coincident peak events. By considering versions \texttt{b} and \texttt{c}, we see that no CP is missed anymore. The main difference between those two versions lies in the segregation of the alerts: version \texttt{b} provides alerts from all four colors, while version \texttt{c} leads to a stronger signal (identifying all CP with a probability above 0.8) at the expense of a higher number of alerts (although restricting only to red ones, they lead to a lower number of alerts). 

\begin{table}[htb!]
    \centering
    \begin{tabular}{||c||c|c||c|c||}
    \hline 
       \multirow{2}{*}{Strategy}  & \multicolumn{2}{c||}{Simple (\texttt{S})}   &   \multicolumn{2}{c||}{Color (\texttt{C})} \\
       \cline{2-5} 
         & \# alerts &  \# CP &  \# alerts (\texttt{R}/\texttt{O}/\texttt{Y}/\texttt{G}) & \# CP (\texttt{R}/\texttt{O}/\texttt{Y}/\texttt{G})\\
       \hline 
        \texttt{1a} & 1.4 & 0.3 & 3.4 (0.7/0.5/0.6/1.6) & 0.6 (0.2/0.1/0.0/0.3) \\
        \texttt{1c} & 6.2 & 0.7 & 10.5 (2.0/2.3/2.6/3.6) & 1.0 (0.5/0.2/0.1/0.2) \\ 
        \texttt{1d} & 13.0 & 1.0 &  13.4 (8.2/3.6/2.7/5.2) & 1.0 (1.0/0.0/0.0/0.0)  \\
        \texttt{2a} & 1.2 & 0.3 & 1.9 (0.5/0.5/0.3/1.6) & 0.6 (0.1/0.2/0.0/0.3) \\ 
        \texttt{2c} & 5.4 & 0.7 & 9.8 (1.7/1.9/2.9/3.3) & 1.0 (0.5/0.2/0.2/0.1) \\ 
        \texttt{2d} & 12.1 & 1.0 & 19.3 (7.8/3.3/2.6/5.6) & 1.0 (1.0/0.0/0.0/0.0) \\
        \texttt{3a} & 1.2 & 0.3 & 3.0 (0.5/0.5/0.3/1.7) & 0.6 (0.1/0.2/0.0/0.3) \\
        \texttt{3c} & 5.4 & 0.7 & 9.7 (1.8/2.0/2.7/3.2) & 1.0 (0.5/0.2/0.1/0.2)\\
        \texttt{3d} & 12.4 & 1.0 & 19.2 (7.2/3.1/3.0/5.4) & 1.0 (1.0/0.0/0.0/0.0) \\
    \hline 
    \end{tabular}
\caption{Average performance statistics of selected strategies for the 1CP program of PSEG. The averages are computed over 10 years, the color strategy shows the overall average performance broken down according to the four colors: red \texttt{R}, orange \texttt{O}, yellow \texttt{Y} and green (\texttt{G})}
\label{tab:all_pseg_perf_cp_day_stats}
\end{table}

\begin{figure}[htb!] 
    \centering
    \includegraphics[width = 16cm]{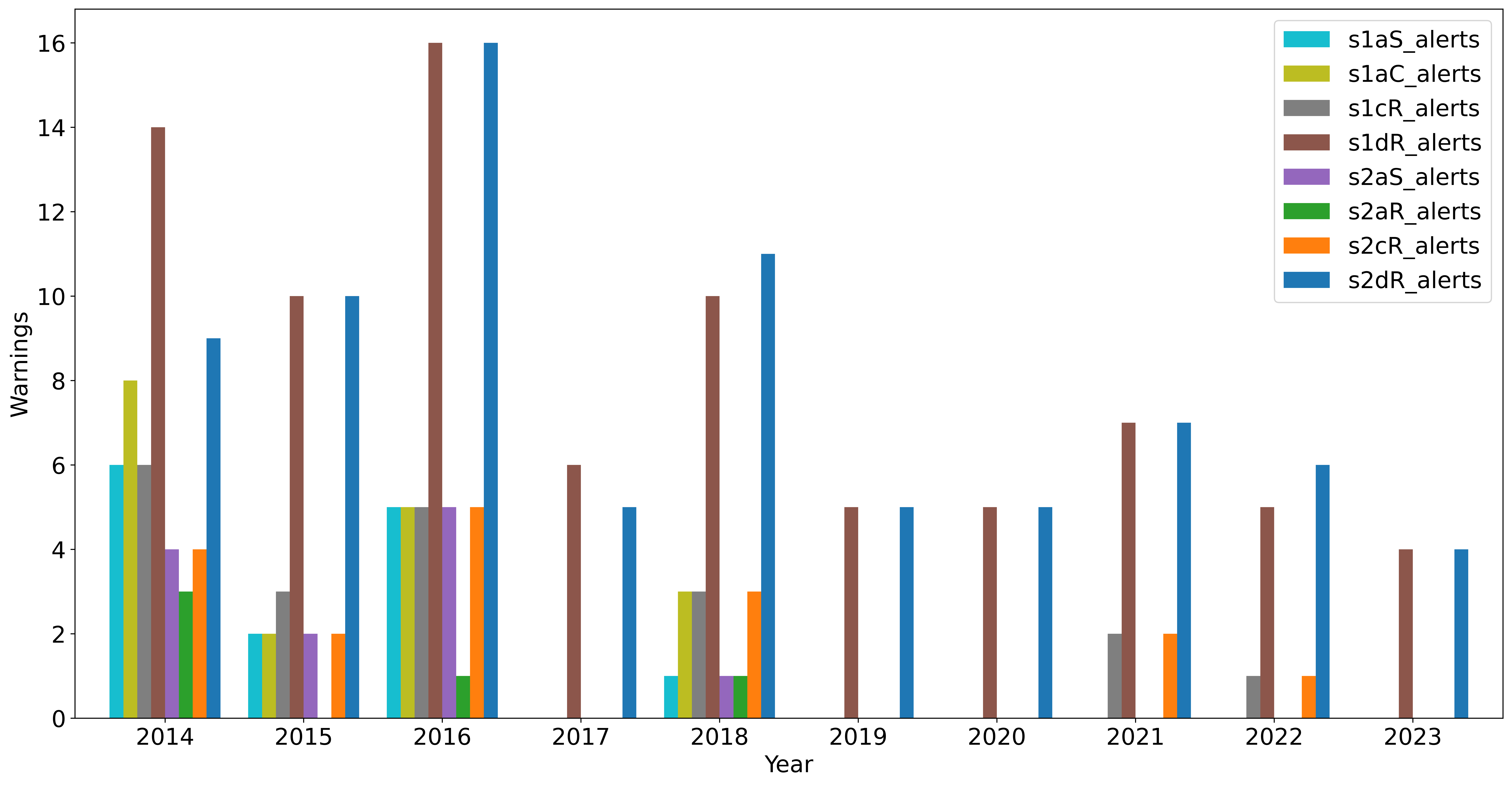}
    \caption{Number of alerts sent for selected strategies broken down per year in the testing period for PSEG 1CP program}
    \label{fig:ps_alert_distribution} 
\end{figure}

Looking at a yearly break-down of the warnings, we see that some of the strategies, e.g. \texttt{1a} or \texttt{2a}, are not sending any warnings in some years including 2017, 2019, and 2020, thus leading to missing the coincident peaks. On the contrary, the strategies \texttt{1cR} and \texttt{2cR} will send out warnings every year and lead to catch the true coincident peaks, as indicated in Table \ref{tab:all_pseg_perf_cp_day_stats}. 

\newpage 

\subsubsection{PJM} Table \ref{tab:pjm_perf_cp_day} shows that the strategies are efficient in the detection of the true CP days with an average of caught CPs of 4.9 irrespective of the strategy. In particular, we notice that most of the 5CPs over the years are caught by the most restrictive signal \texttt{R} defined for probabilities being above $0.8$. In accordance with Subsection \ref{subsec:thresh_5cp}, the average number of warnings does increase when the threshold is lower for a given threshold method. Given the yearly alert breakdown, any method involving no thresholding (\texttt{1}) generates too many false alerts and should therefore be set aside. Preference is to be given to the other strategies that provide at most 15 alerts and catching all the true CPs. 
 
\begin{table}[htb!]
    \centering
    \begin{tabular}{||c||c|c||c|c||}
    \hline 
       \multirow{2}{*}{Strategy}  & \multicolumn{2}{c||}{Simple}  &   \multicolumn{2}{c||}{Color} \\
       \cline{2-5} 
        & \# alerts &  \# CP &  \# alerts (\texttt{R}/\texttt{O}/\texttt{Y}) & \# CP (\texttt{R}/\texttt{O}/\texttt{Y}) \\
       \hline 
        \texttt{1a}  & 26.0 & 4.9 & 26.9 (23.0/2.5/1.4) & 4.9 (4.3/0.6/0.0)  \\
        \texttt{1b} & 30.1 & 5.0 & 31.2 (26.1/3.0/2.1) & 5.0 (5.0/0.0/0.0) \\ 
        \texttt{1c} & 37.7 & 5.0 & 38.7 (33.9/3.0/1.8) & 5.0 (5.0/0.0/0.0) \\
        \texttt{2a} & 12.7 & 4.9 & 13.3 (9.5/2.5/1.3) & 4.9 (4.2/0.7/0.0)  \\
        \texttt{2b} & 15.2 & 5.0 & 16.2 (11.4/2.5/2.3) & 5.0 (4.8/0.2/0.0) \\ 
        \texttt{2c} & 17.8 & 5.0 & 18.2 (13.6/2.9/1.7) & 5.0 (4.8/0.2/0.0) \\
        \texttt{3a} & 14.7 & 4.9 & 16.0 (12.1/2.3/1.6) & 4.9 (4.3/0.6/0.0) \\ 
        \texttt{3b} & 18.9 & 5.0 & 20.1 (15.8/2.0/2.3) & 5.0 (5.0/0.0/0.0) \\ 
        \texttt{3c} & 23.1 & 5.0 & 24.1 (20.1/2.1/1.9) & 5.0 (5.0/0.0/0.0) \\
        \hline 
    \end{tabular}
    \caption{Average performance statistics of selected strategies for the 5CP program of PJM using the forecast (11). The averages are computed over 10 years, the color strategy shows the overall average performance broken down according to the three colors: red \texttt{R}, orange \texttt{O}, and yellow \texttt{Y})}
    \label{tab:pjm_perf_cp_day}
\end{table}

\begin{figure}[htb!]
    \centering
    \includegraphics[width=16cm]{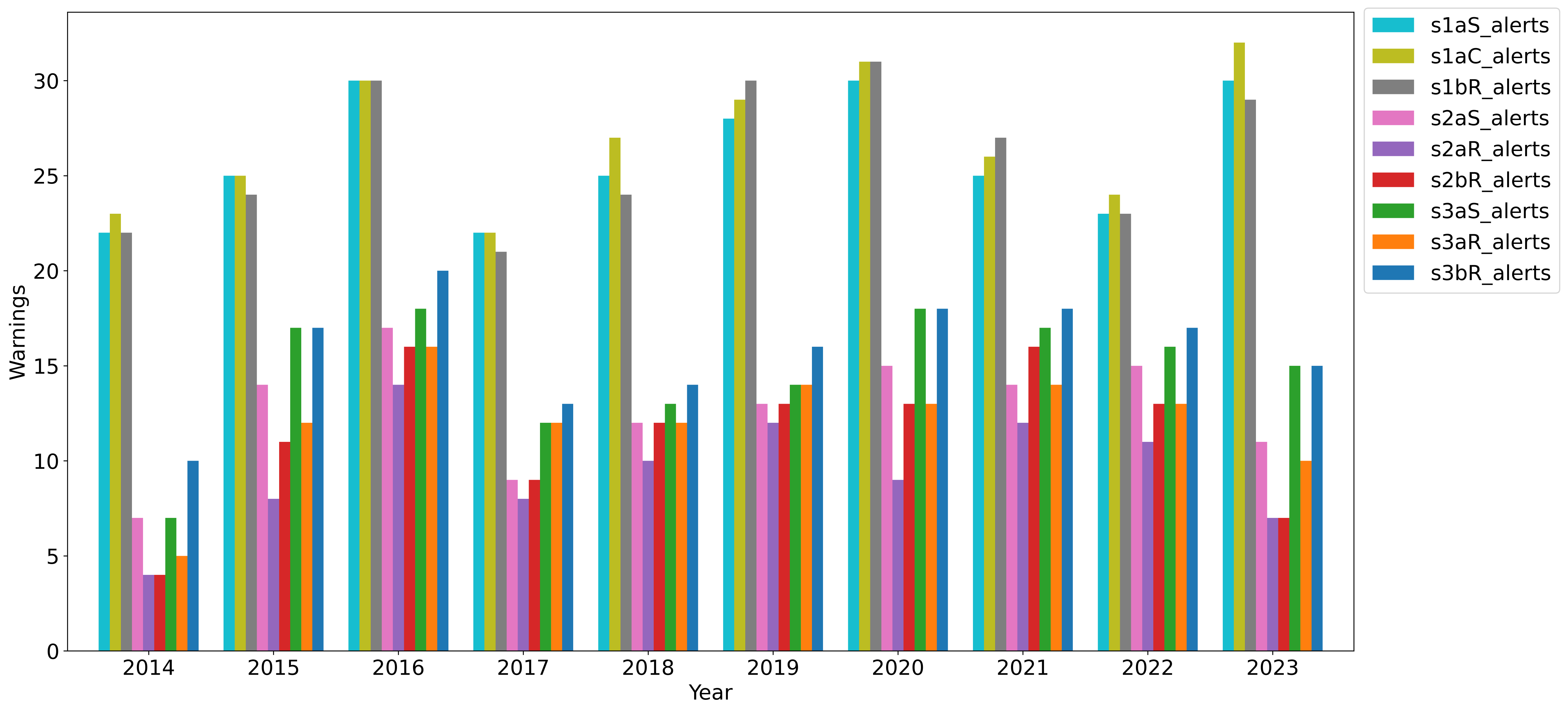}
    \caption{Breakdown per year of the number of warnings generated by selected strategies for PJM 5CP program}
    \label{fig:pjm_perf_cp_day_year}
\end{figure}

\subsubsection{Practical implications}
The parameters to be used in the considered strategies will heavily depend on the consumer profile and preferences. Indeed, the fine-tuning of these parameters will have to take into account several factors. For instance, the number of warnings will need to be determined in function of the impact on productivity or comfort of users under load curtailment. Too many warnings will lead to numerous dispatch that can lead to fatigue of the users. On the other hand, the actual strategy also depends on the risk aversion towards missing a coincident peak event. 

 
\subsection{Performance on Coincident Peak Hour Prediction}

We analyze the performance of our predictor on the CP hour prediction for some of the strategies presented above. To that end, we consider that the estimator perfectly caught the coincident peak hour if it occurs at the same time as the hour with the highest assigned probability and then consider an error $k$ if it has the probability rank $k+1$ for $k \geq 1$. We clip this error at $k=4$ in the sense that if the CP hour is not detected to be among the hours with the top 4 probabilities, it gets assigned the clipped value of 4. 


\subsubsection{NYISO} 

Across strategies, at least 78\% of the alerts assign one of the two highest probabilities to the real daily maximum load hour. When looking at the true coincident peaks, a duration of 2 hours of load curtailing would have led to capturing all the true coincident peaks, and a duration of 1 hour would have led to capturing more than 56\% of the detected CPs. 

\begin{figure}[htb!] 
    \centering
    \begin{subfigure}[t]{0.45\textwidth}
        \includegraphics[width = \textwidth]{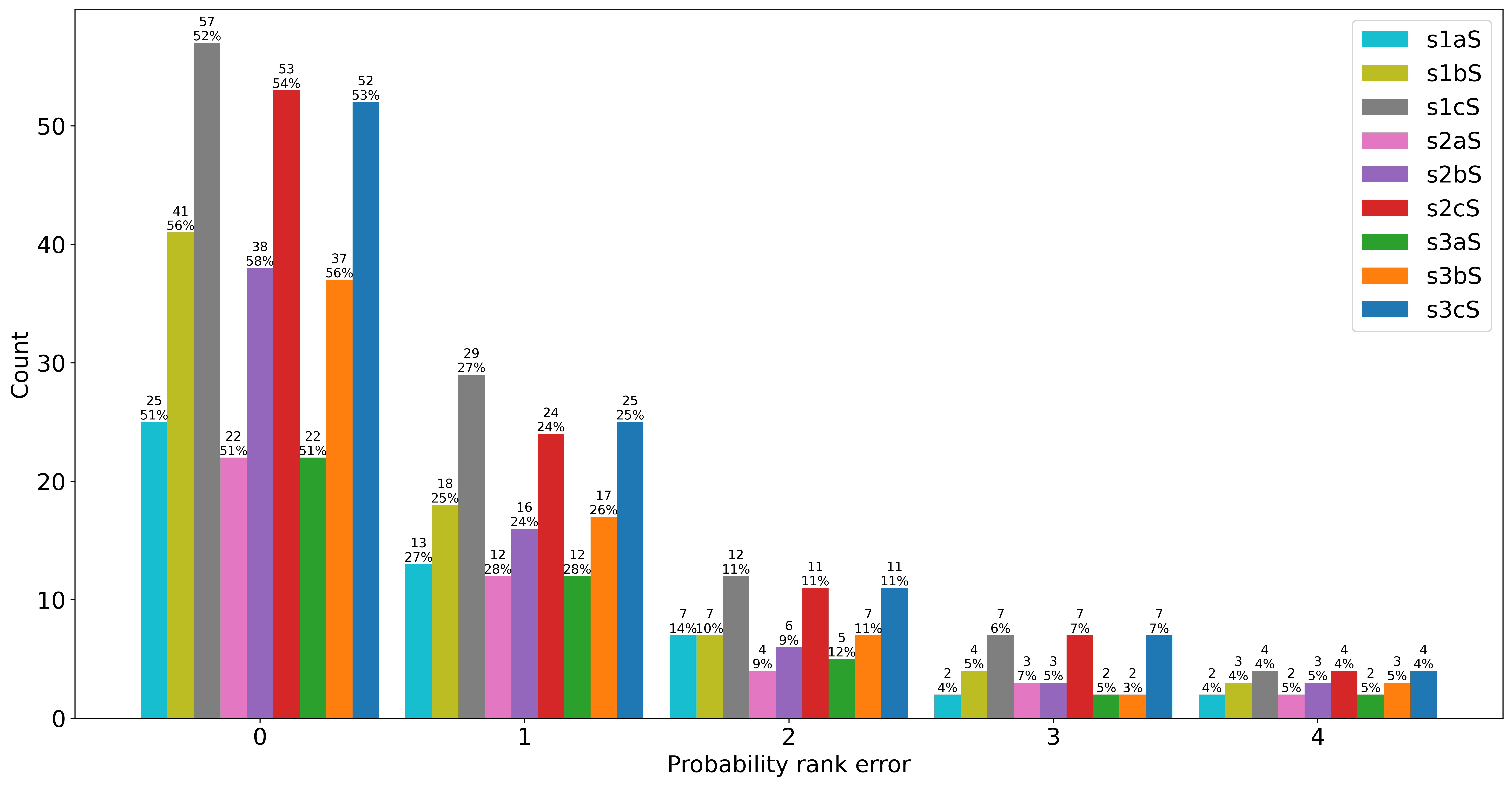}
        \caption{Histogram of rank error for total number of alerts per strategy}
    \end{subfigure}
    ~ 
    \begin{subfigure}[t]{0.45\textwidth}
        \includegraphics[width = \textwidth]{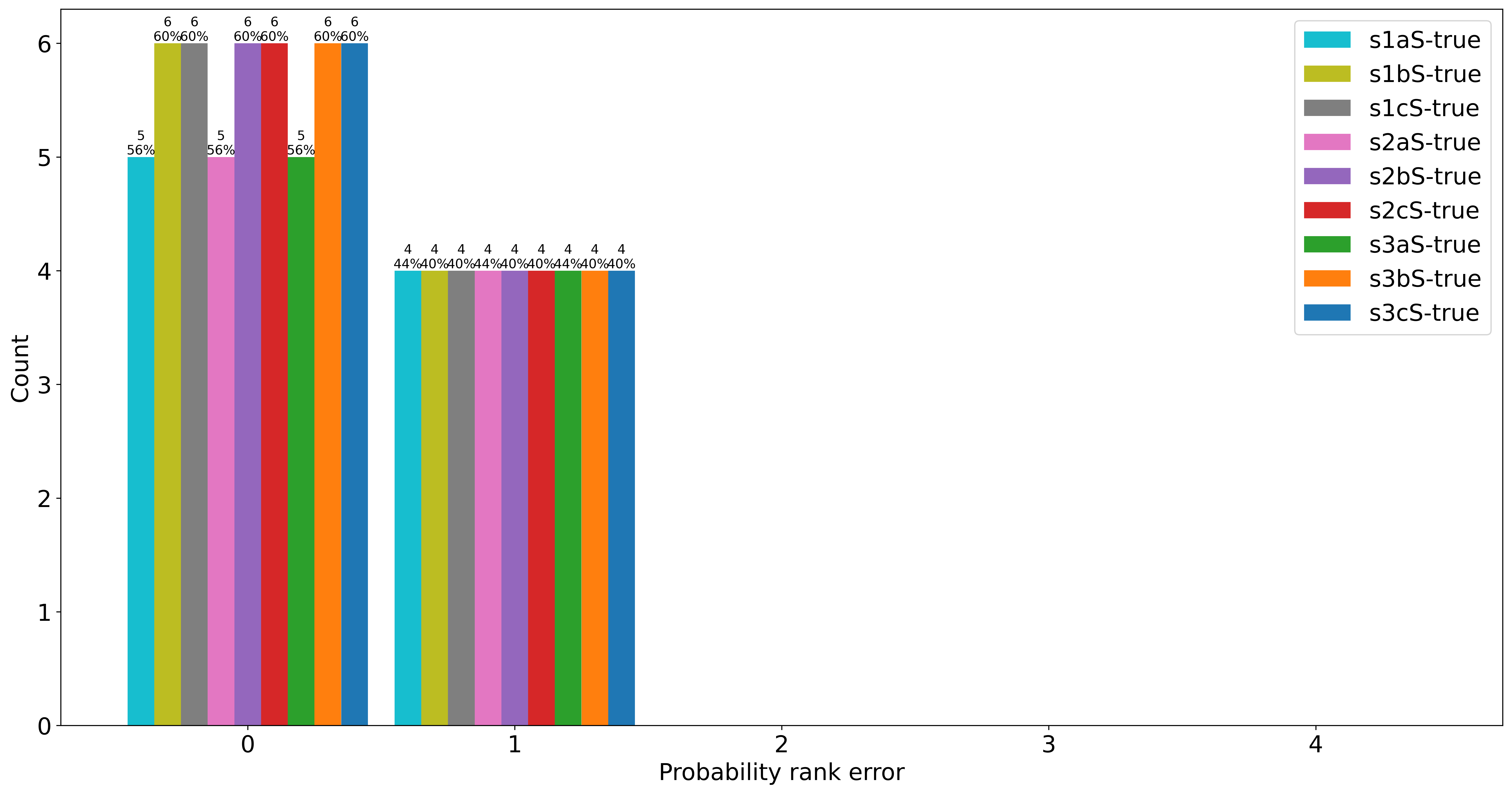}
        \caption{Distribution of rank error for total number of CP caught per strategy}
    \end{subfigure}
    \caption{Histogram and distribution of probability rank error for NYSIO for selected strategies. Each label first shows the number of alerts with the given error and the respective proportion with respect to the total number of alerts generated by each strategy.}
    \label{fig:nyiso_rank_distribution} 
\end{figure}


\subsubsection{PSEG}

\begin{figure}[htb!] 
    \centering
    \begin{subfigure}[t]{0.45\textwidth}
        \includegraphics[width = \textwidth]{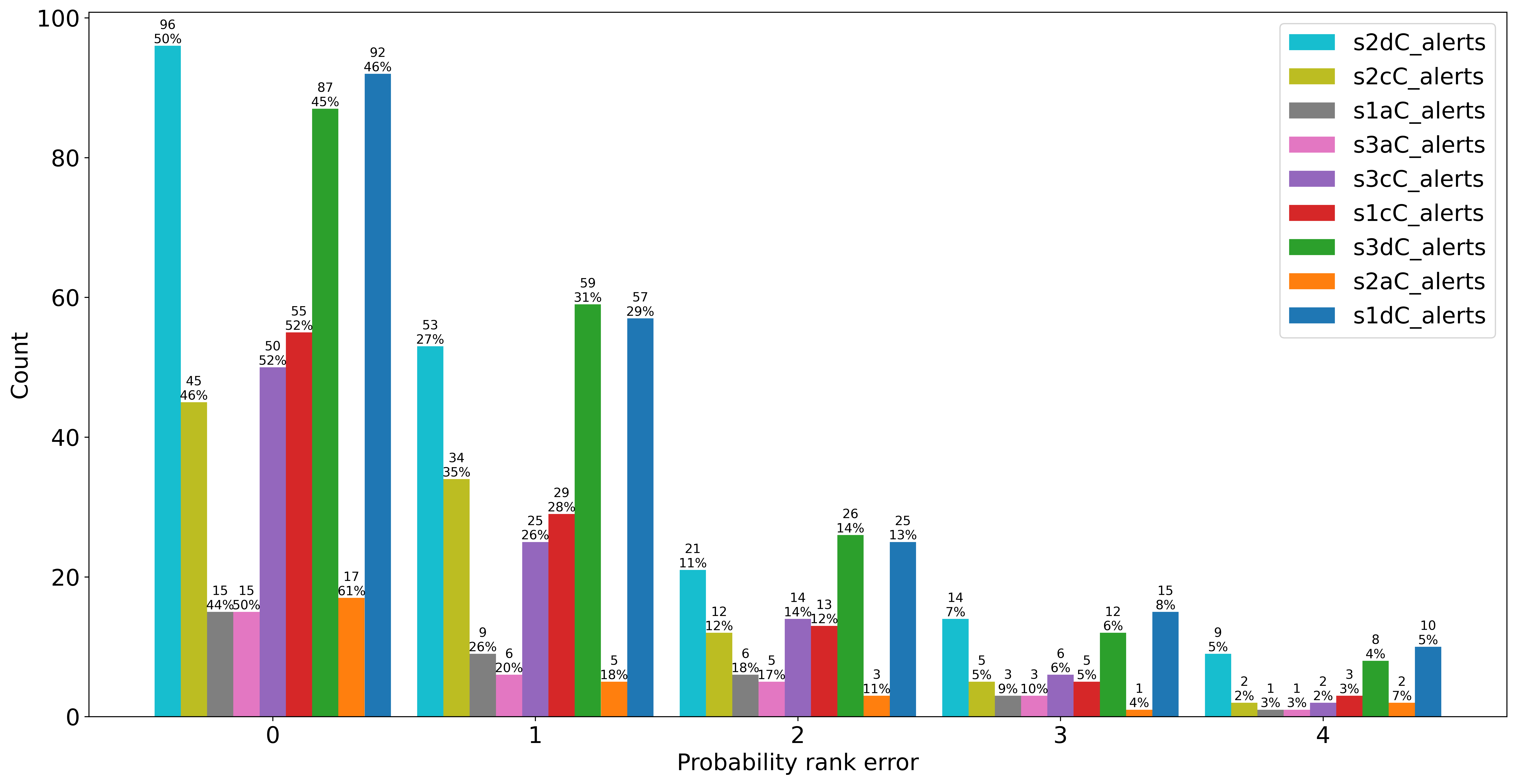}
        \caption{Histogram of rank error for total number of alerts per strategy}
    \end{subfigure}
    ~ 
    \begin{subfigure}[t]{0.45\textwidth}
        \includegraphics[width = \textwidth]{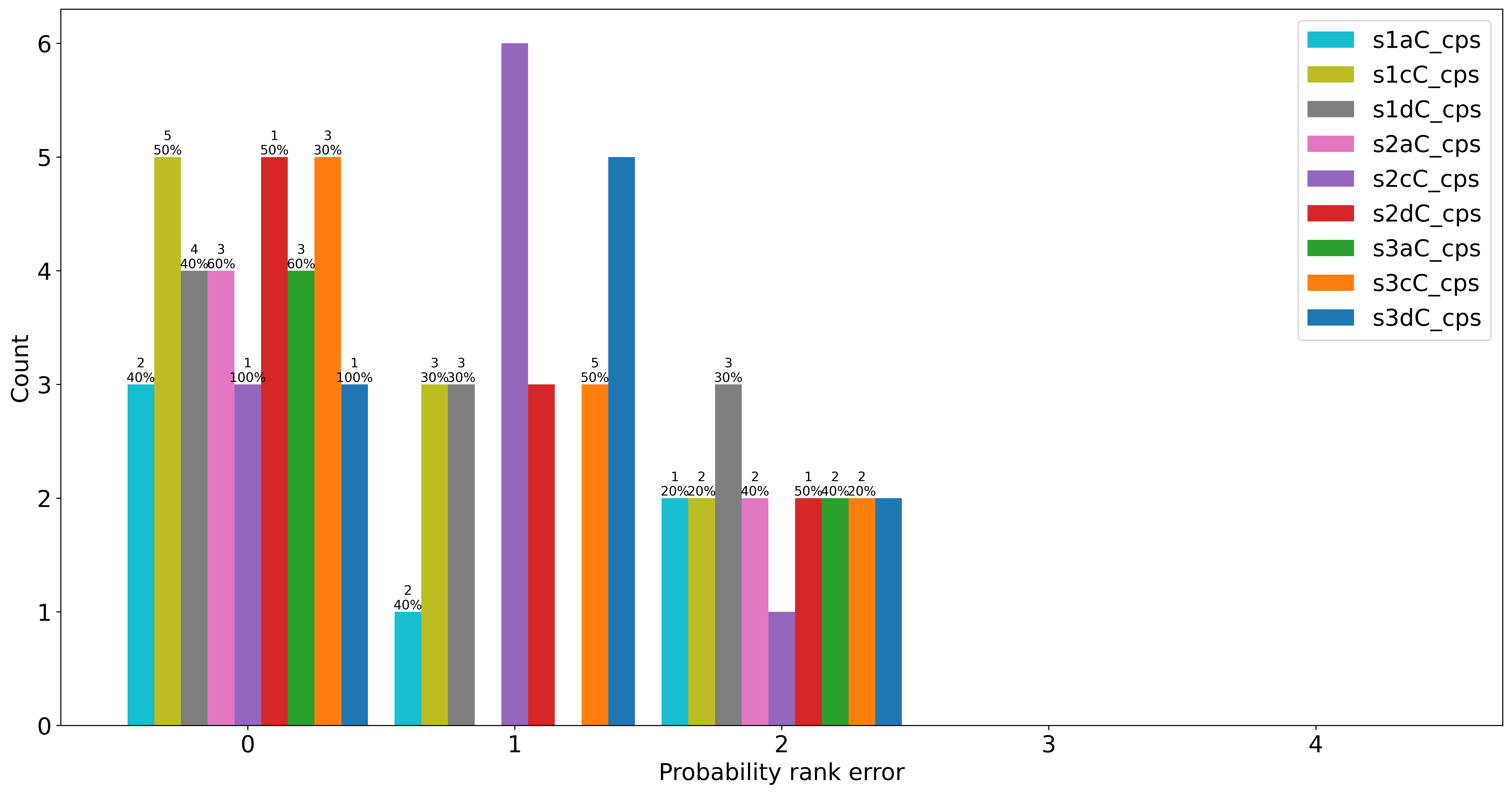}
        \caption{Distribution of rank error for total number of CP caught per strategy}
    \end{subfigure}
    \caption{Histogram and distribution of probability rank error for PSEG for selected strategies.} 
    \label{fig:pseg_rank_distribution} 
\end{figure}

We notice that at least 70\% of all generated alerts end up assigning one of the two highest probabilities to the true daily maximum and at least 88\% assigning one of the three highest probabilities. Looking only at the caught coincident peak events, it shows that the CP hours are all assigned one of the top three probabilities. 

\subsubsection{PJM}

Across strategies, more than 80\% of the alerts assign one of the two highest probabilities to the real daily maximum load hour. When looking at the caught final coincident peaks, a duration of 3 hours of load curtailing would have led to capturing at least 94 \% of the the detected final coincident peaks, and a duration of 2 hours would have led to capturing 84\% of them. The study shows that a load curtailment duration of 4 hours would enable to reduce the load on all detected coincident peaks, which are between 48 and 50 of the real CPs depending on the adopted strategy. Overall, the detection of the daily maximum load hour on the final coincident peaks is performing really well with no true CP hours being assigned a probability that is not among the top 4 probabilities. 

\begin{figure}[htb!] 
    \centering
    \begin{subfigure}[t]{0.45\textwidth}
        \includegraphics[width = \textwidth]{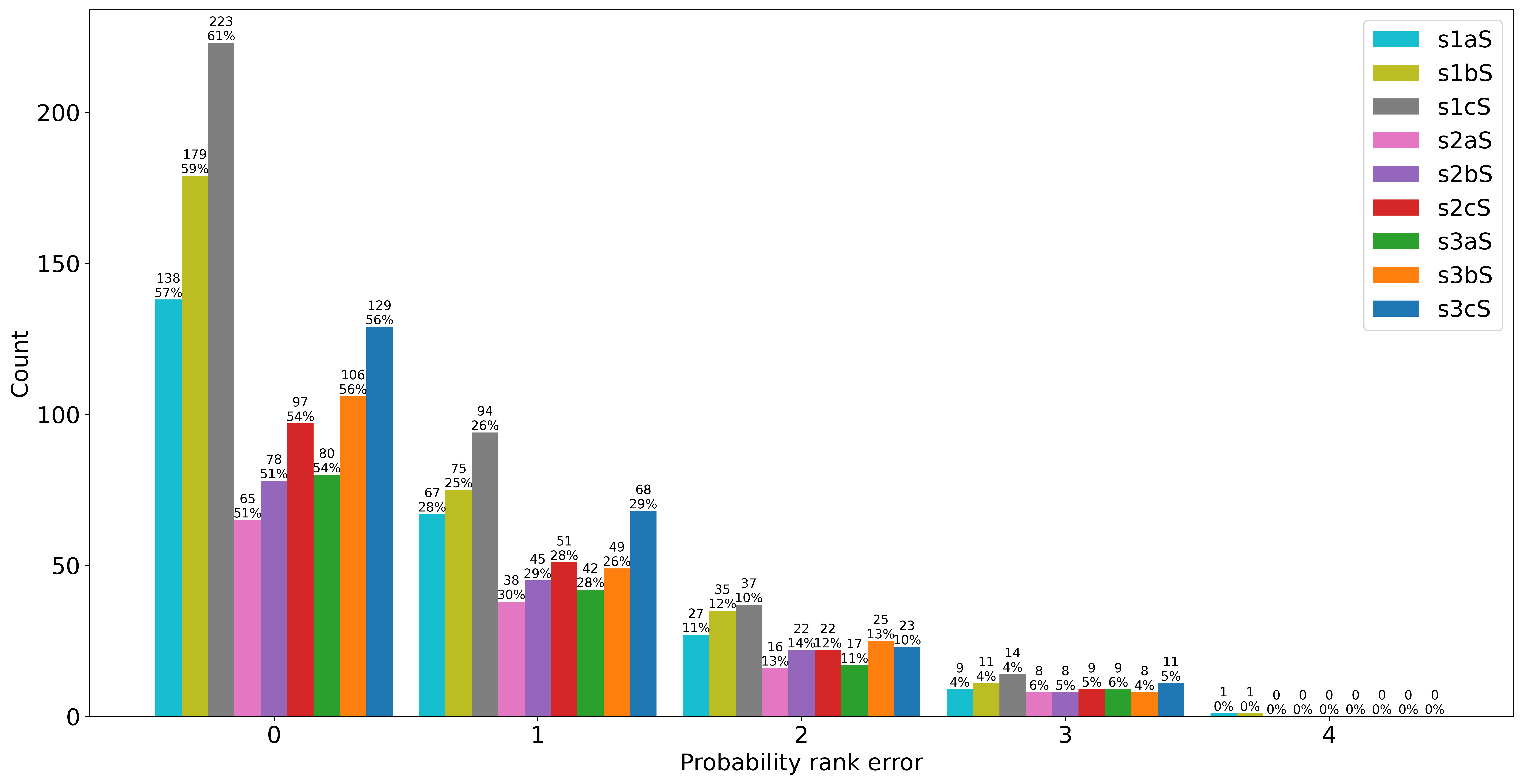}
        \caption{Histogram of rank error for total number of alerts per strategy}
    \end{subfigure}
    ~ 
    \begin{subfigure}[t]{0.45\textwidth}
        \includegraphics[width = \textwidth]{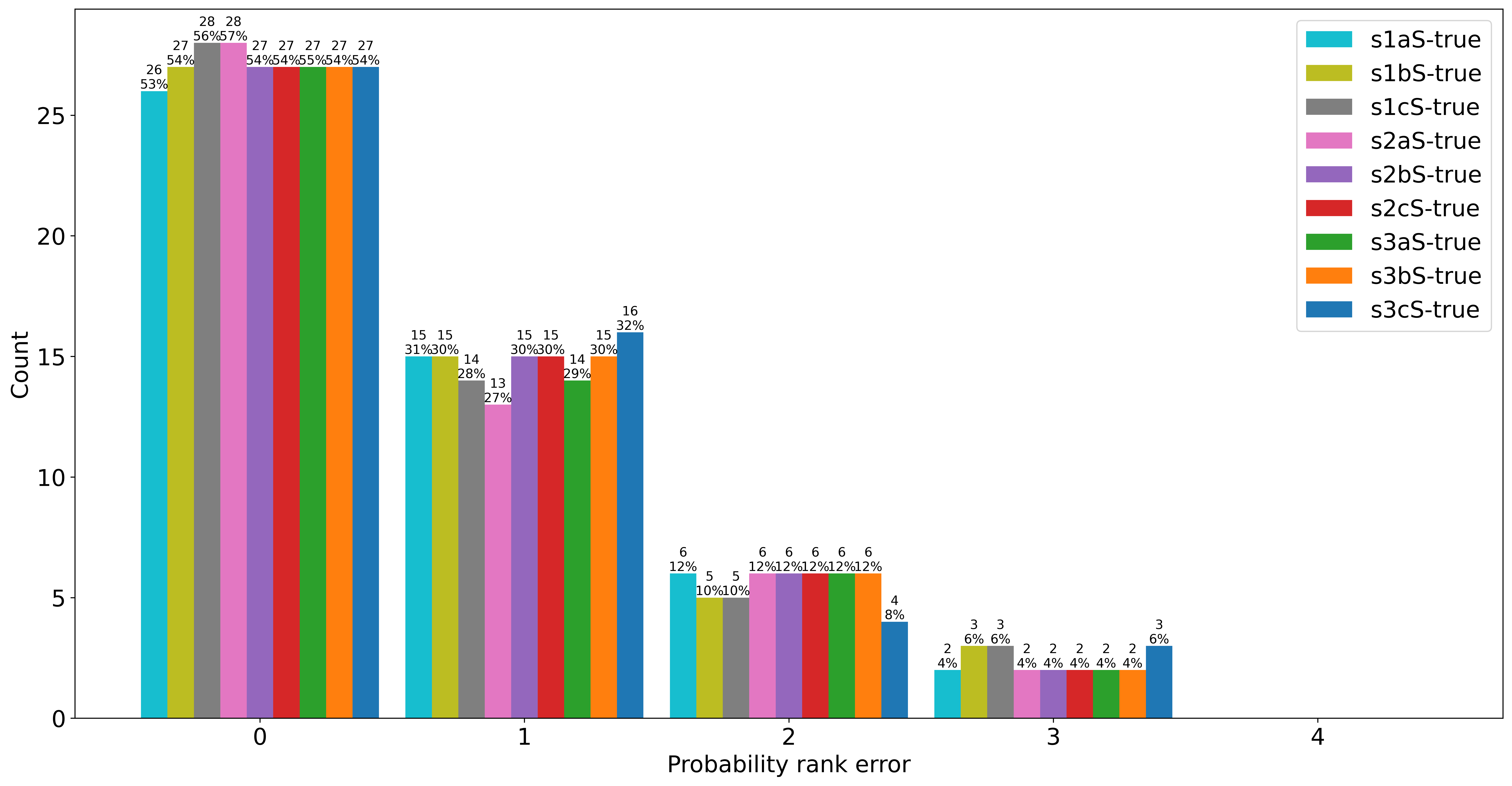}
        \caption{Distribution of rank error for total number of CP caught per strategy}
    \end{subfigure}
    \caption{Histogram and distribution of probability rank error for PJM for selected strategies. }
    \label{fig:pjm_rank_distribution} 
\end{figure} 

\subsubsection{Practical implications}

The load curtailment effort will need to last a certain duration on alert days in order to capture the coincident peak: the larger the time window around the identified maximum hour is, the more certain the load reduction on the actual coincident peak hour is. In particular, the above analysis can help consumers to size any local generator or energy storage solution, as well as quantify the potential economic benefits and opportunity costs. 


\section{Conclusion} \label{sec:conclu} 

This work presents a scenario generation engine capable of handling two different cases, depending on the availability of an electric load forecast. The first application includes the coincident peak program computed for a single zone or region, e.g. an ISO, for which forecasts of electric load are publicly available. Our estimators are successful in predicting the day and hour of the events for both the single and the multiple coincident peak settings. Moreover, our back-tests also demonstrate the relative generality and transportability of our method.  The second application concerns the case of a CP progra in a zone for which there is no available forecast, forcing us to rely on a proxy, for example the load  forecasts for a different region. To overcome this hurdle, we develop a conditional scenario generation engine based on a joint stochastic model that can be fitted to the available data. Finally, we develop a sequential scheme with adaptive thresholding that estimates the probability of coincident peak events occurring over the next day and at which hours. We analyze the performance of several such strategies and draw some practical implications for load curtailment dispatch time window and sizing characteristics of BESS technologies. 

\printbibliography

\begin{appendix}

\section{Description of PJM Network Integration Transmission Service (NITS) Rates by utility company} \label{app:PJM_NITS}

\begin{table}[htb!]
    \centering
    \resizebox{0.9\columnwidth}{!}{%
    \begin{tabular}{|c|c|c|c|c|}
    \hline
       \multicolumn{2}{|c|}{\textbf{Transmission Zone}} & \textbf{Region} & \multicolumn{2}{|c|}{\textbf{CP Methodology}}\\ \hline \hline 
       \textbf{Utility Company} & \textbf{Acronym} & \textbf{Mid-Atlantic}  & \textbf{\#CP (System)} & \textbf{Period} \\ \hline \hline 
       Allegheny Power	& APS/AP/FE South & \xmark & 1CP  (LSE) & Jun 1 -- Sep 30 \\ 
        \hline 
        American Electric Power	& AEP & \xmark & 1CP  (LSE)  & Nov 1 -- Oct 31 \\ 
        \hline 
       Atlantic City Electric & AECO/AE & \cmark &  5CP (LSE)  & Jun 1 -- Sep 30 \\ 
             \hline 
         \multirow{2}{*}{American Transmission Service Inc} & \multirow{2}{*}{FE/ATSI} & \multirow{2}{*}{\cmark} &  \multirow{2}{*}{ 5CP (LSE) } & Jun 1 -- Sep 30 \\ 
       & & & &  Dec 1 -- Mar 31 \\
        \hline 
       Baltimore Gas and Electric & BGE/BC & \cmark & 5CP (PJM) & Jun 1 -- Sep 30 \\
              \hline 
        Commonwealth Edison & ComEd & \xmark &  1CP (LSE) & Jun 1 -- Sep 30 \\
              \hline 
        Delmarva Power and Light Company & &\xmark & 5CP (LSE) & Nov 1 -- Oct 31 \\
              \hline 
        Dominion Virginia Power & DOM/DVP & \xmark & 12CP (LSE)& Nov 1 -- Oct 31 \\ 
              \hline 
        Duke Energy Ohio and Kentucky & DEOK & \cmark &  1CP (LSE)& Nov 1 -- Oct 31\\ 
              \hline 
        Duquesne Light & DUQ/DLCO & \xmark &  1CP (LSE)& Jan 1 -- Dec 31 \\ 
              \hline 
         East Kentucky Power Cooperative &  EKPC &  \xmark & 1CP (PJM)  &  Jun 1 -- Sep 30   \\        \hline 

       \multirow{2}{*}{Jersey Central Power and Light} & \multirow{2}{*}{JCPL/JC} & \multirow{2}{*}{\cmark} &  \multirow{2}{*}{5CP (LSE)} & Jun 1 -- Sep 30 \\ 
       & & & &  Dec 1 -- Mar 31 \\
              \hline 
        \multirow{2}{*}{Metropolitan Edison} & \multirow{2}{*}{MetEd/ME} & \multirow{2}{*}{\cmark} &  \multirow{2}{*}{5CP (LSE)} & Jun 1 -- Sep 30 \\ 
       & & & & Dec 1 -- Mar 31 \\
              \hline 
        Old Dominion Electric Cooperative & ODEC & \xmark &  5CP (LSE) & Nov 1 -- Oct 31 \\       
        \hline 
       PECO Energy & PECO/PE & \cmark &  1CP (LSE) & Nov 1 -- Oct 31 \\ 
              \hline 
        \multirow{2}{*}{Pennsylvania Electric Company} & \multirow{2}{*}{PENELEC/PN} & \multirow{2}{*}{\cmark} &  \multirow{2}{*}{5CP (LSE)} &  Jun 1 -- Sep 30  \\
       & & & & Dec 1 -- Mar 31 \\
        \hline 
       PPL Electric Utilities & PPL/PL & \cmark & 5CP (LSE) & Nov 1 -- Oct 31   \\ 
              \hline 
        Potomac Electric Power &  Potomac/PEPCO& \cmark &   1CP (LSE) & Jun 1 -- Sep 30  \\ 
       \hline 
       Public Service Electric and Gas & PSEG/PS & \cmark & 1CP (LSE) & Jun 1 -- Sep 30  \\ 
       \hline 
    \end{tabular}
    }
    \caption{Cost Allocation Methodology for some PJM Utility companies. The methodology is indicated by the column \texttt{\#CP} where the coincident peaks are computed from the load of the system in parenthesis over the time interval indicated by the column \texttt{Period}. LSE denotes the Load Serving Entity or utility and PJM the RTO. }
    \label{tab:PJM_LSE_CP}
\end{table}

\end{appendix}

\end{document}